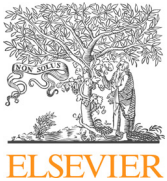

Contents lists available at ScienceDirect

# Neurocomputing

journal homepage: www.elsevier.com/locate/neucom

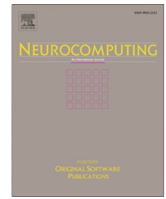

# Stability analysis of nontrivial stationary solution and constant equilibrium point of reaction–diffusion neural networks with time delays under Dirichlet zero boundary value

Ruofeng Rao [a,1], Jialin Huang [b], Xiaodi Li [c]

[a] Department of Mathematics, Chengdu Normal University, Chengdu 611130, China
[b] Department of Mathematics, Sichuan Sanhe Vocational College, 646200 Luzhou, Sichuan, China
[c] School of Mathematics and Statistics, Shandong Normal University, Jinan 250014, China

A B S T R A C T

In this paper, Lyapunov–Razumikhin technique, design of state-dependent switching laws, a fixed point theorem and variational methods are employed to derive the existence and the unique existence results of globally exponentially stable (positive) stationary solution of delayed reaction–diffusion cell neural networks under Dirichlet zero boundary value, including the global stability criteria *in the classical meaning*. Next, sufficient conditions are proposed to guarantee the global stability invariance of ordinary differential systems under the influence of diffusions. New theorems show that the diffusion is a double-edged sword in judging the stability of diffusion systems. Besides, an example is constructed to illuminate that any non-zero constant equilibrium point must be not in the phase plane of dynamic system under Dirichlet zero boundary value, or it must lead to a contradiction. Next, under Lipschitz assumptions on active function, another example is designed to prove that the small diffusion effect will cause the essential change of the phase plane structure of the dynamic behavior of the delayed neural networks via a Saddle point theorem. Finally, a numerical example illustrates the feasibility of the proposed methods. It is worth mentioning that some interesting mathematical problems are **originally** put forward in the last chapter.

© 2021 The Authors. Published by Elsevier B.V. This is an open access article under the CC BY-NC-ND license (http://creativecommons.org/licenses/by-nc-nd/4.0/).

## 1. Introduction

Firstly, we recalled the reason why we need to study the stability of reaction–diffusion neural networks system. In 1988, inspired by cellular automata, Chua and Yang proposed a new neural network based on Hopfield network, i.e. cellular neural network (CNN), which is formed by a number of cells with the same structure after a well-organized combination [22,23]. Each neuron in the network will automatically choose to connect with the nearest neuron. Because of its local connectivity, CNN is especially suitable for ultra large scale integrated circuit implementation. The characteristics of the above cellular neural network make it widely used in pattern recognition, image processing, signal processing and other fields. The main function of cellular neural network is to transform an input image into a corresponding output image. For example, the existing target motion direction detection, edge detection, and connected slice detection all use this function. In order to achieve these functions, the cellular neural network must be completely stable, that is, all output trajectories must converge to a stable equilibrium point. So the stability of cellular neural network has become a hot topic [24–26]. As we all know, time delay may destroy the stability of the system and lead to oscillation, bifurcation, chaos and other phenomena, thus changing the characteristics of the system. In cellular neural networks, time delay is inevitable. For example, there are cell delay, transmission delay and synapse delay in biological neural network [27]. As pointed out in [28] that many pattern formation and wave propagation phenomena that appear in nature can be described by systems of coupled nonlinear differential equations, generally known as reaction–diffusion equations. These wave propagation phenomena are exhibited by systems belonging to very different scientific disciplines. Besides, the interactions arising from the space-distributed structure of the multilayer cellular neural networks can be seen as diffusion phenomenon [3,29,34–36]. As pointed out by Qiankun Song and Zidong Wang in [34], many pattern formation and wave propagation phenomena that appear in nature can be described by systems of coupled nonlinear differential equations, generally known as reaction–diffusion equations. These wave propagation phenomena are exhibited by systems belonging to very different scientific disciplines. The reaction–diffusion effects, therefore, cannot be neglected in both biological and man-made neural networks, especially when electrons are moving in non-even electromagnetic field. Moreover, although the diffusion coefficients may be very small, the topological structure of the phase plane of the dynamic behavior of the following reaction–diffusion system (1.1) is likely to change substantially from a constant equilibrium point of the following system (1.3) to multiple stationary solutions of the reaction–diffusion system (1.1). Therefore, many global stability results of delayed neural networks in the form of ordinary differential equations may only be locally asymptotical stability criteria in real engineering. *Unfortunately, such an example has not been constructed for the time being. But that doesnt mean there are no such examples, which may become an open problem hereafter. On the other hand, fortunately, this paper has proposed the conditions guaranteeing the global stability invariance of ordinary differential systems under the influence of diffusions in the meaning of* Definition 1 *(see Corollary 3.4).*

---
[1] The work was jointly supported by the Application basic research project of science and Technology Department of Sichuan Province (2020YJ0434) and the Major Scientific research projects of Chengdu Normal University in 2019 (CS19ZDZ01).
   E-mail addresses: ruofengrao@163.com (R. Rao), jialink2880@163.com (J. Huang), lxd@sdnu.edu.cn (X. Li).

Next, we shall point out the fact that the stability results in previous literature involved to delayed reaction–diffusion neural networks make it unnecessary to study the reaction–diffusion system (partial differential equations model), but only its corresponding ordinary differential equations model. What's the problem?

For a long time, the stability of the reaction diffusion neural networks was investigated in many literatures [1–5], in which the stability of the constant equilibrium point was studied. For example, in [1], the following cellular neural networks with time-varying delays and reaction–diffusion terms was considered,

$$\frac{\partial y(t,x)}{\partial t} = \sum_{q=1}^{m} \frac{\partial y(t,x)}{\partial x_q} \left( D_q \frac{\partial y(t,x)}{\partial x_q} \right)$$
$$- Cy(t,x) + Ag(y(t,x)) + Bg(y(t-\tau(t),x)) + J, (t,x) \in \mathbb{R}_+ \times \Omega, \tag{1.1}$$

Next, the authors of [1] defined the equilibrium point of the time-delayed reaction–diffusion system (1.1) as the constant vector $y^*$ satisfying

$$Cy^* = Ag(y^*) + Bg(y^*) + J. \tag{1.2}$$

Here, we have to say, the equilibrium point $y^*$ is also the equilibrium point of the following ordinary differential equations corresponding to the time-delayed partial differential Eqs. (1.1),

$$\frac{dx(t)}{dt} = -Cx(t) + Ag(x(t)) + Bg(x(t-\tau(t))) + J, t \in \mathbb{R}_+, \tag{1.3}$$

Due to the Poincare inequality, we see, the diffusion items actually promote the stability of the reaction diffusion system (1.1). That is, we only need to study the ordinary differential Eqs. (1.3) because the stability criteria of the ordinary differential Eqs. (1.3) must make the system (1.1) stable. In other words, the stability of the reaction–diffusion model does not need to be studied because it is included in the stability of its corresponding ordinary differential equations model.

So we need to ask where the problem is? *The answer lies in the fact pointed out in this paper. In Theorem 3.3 and Corollary 3.4, the constant equilibrium point $u^*$ may become another $u^*(x)$, where $u^*(x) \neq u^*$ in common cases under Dirichlet zero boundary value. Only the zero solution $u^* = 0$ might become one of the stationary solutions of reaction–diffusion system. And so the uniqueness existence conditions of the stationary solution $u^*(x)$ (see Corollary 3.4) illuminate the inconvenient and difficulties due to the inevitable diffusions in practical engineering, which make it more difficult to judge the stability than ordinary differential system.* Usually, the positive stationary solution $u^*(x) > 0$ implies more realistic meanings in neural networks and other systems, such as the financial systems [8,30]. And so Theorem 3.1 and Theorem 3.2 have proposed the existence and the unique existence of the (globally) exponentially stable positive stationary solution in this paper.

On the other hand, Neural Networks usually show the special characteristics of network mode switching, which results in the so-called switched Neural Networks. That is to say, switched Neural Networks are a kind of Neural Networks whose parameters are operated by a switching signal. The authors of [7] designed an state-dependant switching (SDS) law for passivity of switched neural networks subject to stochastic disturbances and time-varying delays, where the proposed SDS law involves both the current state and the previous value of the switching signal. In recent years, switched Neural Networks have received much attention due to their widespread applications in mechanical systems, circuits and power systems, the automotive industry, air traffic control, and many other fields [16,37,38]. In order to describe the switching phenomenon in Neural Networks, the switched Neural Networks have been studied in recent years (see [37,38] and its references therein). As pointed out in [40], the connections between neurons of such networks are not always constant in practice. This usually leads to link failures and new link creation. This sudden change leads to the transition between some different topologies of the network structure and the instability of the Neural Networks. For example, a hot topic on the stability of Neural Networks with random and deterministic switching laws is studied in [40,42].

Motivated by some methods of [1–45], particularly [13–15,18,42], we investigate the stability of the nontrivial stationary solution of switched reaction–diffusion neural networks with time delays. This paper has the following innovations:

★ It is the first paper to study and obtain the existence theorem and unique existence theorem (see Theorem 3.1–3.2, Corollary 3.6–3.7) of globally asymptotically stable nontrivial stationary solution of reaction–diffusion neural networks with time delays under Dirichlet zero boundary value via the comprehensive applications of Lyapunov–Razumikhin technique, design of state-dependent switching laws, a fixed point theorem, variational methods, and construction of compact operators on a convex set. Moreover, such new theorems illuminate originally that the diffusion phenomena is the double-edged sword in judging the stability of delayed reaction–diffusion systems.

★ It is the first paper to propose the conditions guaranteeing the global stability invariance of delayed ordinary differential systems under the influence of diffusions in the meaning of Definition 1 (see Theorem 3.3 and Corollary 3.4). Besides, the conception of the so-called global stability invariance is **originally** proposed in this paper, which leads to the relevant mathematical problem proposed **originally** (see, Problem 4).

★ It is the first paper to design the contradiction results to show that any non-zero constant function must not be a solution of any reaction diffusion neural networks under Dirichlet zero boundary value.

★ It is also the first time to study how the tiny diffusion causes the essential change of the phase plane structure of the dynamic behavior of the delayed neural networks under Lipschitz assumptions on activate functions or signal functions (see Statement 2), which leads to the relevant mathematical problem proposed **originally** (see, Problem 1).

**Notations:**

Throughout this paper, we denote by $I$ the identity matrix with an appropriate dimension, and assume that $D_\sigma$ is a positive definite diagonal matrix for any given $\sigma \in \{1, 2, \ldots, N\}$. Besides, $diag(\cdots)$ stands for a diagonal matrix. $\mathcal{A} > 0$ ($\mathcal{A} \geqslant 0$) means that $\mathcal{A}$ is a real symmetric positive (semi-positive) definite matrix. And denote $|\mathcal{A}| = (|a_{ij}|)_{n \times n}$ for any matrix $\mathcal{A} = (a_{ij})_{n \times n}$. In addition, for any $v = (v_1, v_2, \cdots, v_n)^T, u = (u_1, u_2, \ldots, u_n)^T \in \mathbb{R}^n$, we denote $|v| = (|v_1|, |v_2|, \cdots |v_n|)^T$, and $v \leqslant u$ means that $v_i \leqslant u_i$ for all $i = 1, 2, \ldots, n$. Denote by $\Omega_\sigma$ an open bounded domain in $R^n$ with the smooth boundary $\partial \Omega_\sigma$ for any given $\sigma \in \{1, 2, \ldots, N\}$. For convenience, we denote by $\lambda_{\sigma 1} > 0$ the first positive eigenvalue of the Laplace operator $-\Delta$ on the Sobolev space $W_0^{1,2}(\Omega_\sigma)$.

## 2. System descriptions

Consider the following switched neural networks with time-varying delays and reaction–diffusion terms

$$
\begin{cases}
\frac{\partial y(t,x)}{\partial t} = D_\sigma \Delta y(t,x) - C_\sigma y(t,x) + A_\sigma g(y(t,x)) + B_\sigma g(y(t-\tau(t),x)) + J_\sigma, & (t,x) \in \mathbb{R}_+ \times \Omega_\sigma, \\
y_i(t,x) = 0, t \geqslant 0, x \in \partial \Omega_\sigma, i = 1, 2, \ldots, n,
\end{cases}
\tag{2.1}
$$

where $\Omega_\sigma \subset \mathbb{R}^n$ is a bounded domain with smooth boundary $\partial\Omega_\sigma$, the state variable $y(t,x) = (y_1(t,x), y_2(t,x), \cdots, y_n(t,x))^T$ with $y_i(t,x)$ representing being the state variable of the $i$th neuron in time $t$ and space variable $x$. $J_\sigma = (J_{\sigma 1}, \ldots, J_{\sigma n})^T \in \mathbb{R}^n$ is the constant external input vector, and both $D_\sigma$ and $C_\sigma$ are positive definite diagonal matrices, in which $D_\sigma$ represents the diffusion coefficient matrix, and $C_\sigma$ represents the connection weight matrix of neural network. Besides, $A_\sigma$ and $B_\sigma$ both are the connection weight matrices of neural network. For each $x \in \Omega_\sigma, g(y(t,x)) = (g_1(y_1(t,x)), \ldots, g_n(y_n(t,x)))^T$ represents a time-dependent signal function vector. $\tau(t)$ represents the time delay required for signal transmission from neuron $j$ to neuron $i$, satisfying $0 \leqslant \tau(t) \leqslant \tau$. Assumed that $y^\sigma(x) = (y_1^\sigma(x), \cdots, y_n^\sigma(x))^T$ is a nontrivial stationary solution of reaction–diffusion switched system (2.1), then $y^\sigma(x)$ satisfies two equations of the system (2.1), in addition,

$$
- C_\sigma y^\sigma(x) + A_\sigma g(y^\sigma(x)) + B_\sigma g(y^\sigma(x)) + J_\sigma \neg \equiv 0, \quad x \in \Omega_\sigma.
\tag{2.2}
$$

Of course, the sufficient condition should be given to ensure the existence of such nontrivial stationary solution.

The so-called stationary solution $y^\sigma(x)$ of the system (2.1) in the Sobolev space $W_0^{1,2}(\Omega)$ is that for any $\theta(x) \in W_0^{1,2}(\Omega)$,

$$
D_\sigma \int_\Omega \nabla y^\sigma(x) \nabla \theta(x) dx + C_\sigma \int_\Omega y^\sigma(x) \theta(x) dx - (A_\sigma + B_\sigma)
$$
$$
\int_\Omega g(y^\sigma(x)) \theta(x) dx - \int_\Omega J_\sigma \theta(x) dx = 0,
$$

where the Sobolev space $W_0^{1,2}(\Omega)$ is one of Hilbert spaces.

Set $u(t,x) = y(t,x) - y^\sigma(x)$, then the system (2.1) is translated into the following system:

$$
\begin{cases}
\frac{\partial u(t,x)}{\partial t} = D_\sigma \Delta u(t,x) - C_\sigma u(t,x) + A_\sigma f(u(t,x)) + B_\sigma f(u(t-\tau(t),x)), & (t,x) \in \mathbb{R}_+ \times \Omega_\sigma, \\
u_i(t,x) = 0, t \geqslant 0, x \in \partial \Omega_\sigma, i = 1, 2, \cdots, n,
\end{cases}
\tag{2.3}
$$

where $f(u(t,x)) = g(y(t,x)) - g(y^\sigma(x)), f(u(t-\tau(t),x)) = g(y(t-\tau(t),x)) - g(y^\sigma(x))$. Here, the nontrivial stationary solution $y^\sigma(x)$ of the system (2.1) corresponds to the null solution of the system (2.3).

Besides, we may equip the system (2.3) with the initial value:

$$
u_i(s,x) = \phi_i(s,x), -\tau \leqslant s \leqslant 0, x \in \Omega_\sigma,
\tag{2.4}
$$

where $(\phi_1(s,x), \phi_2(s,x), \cdots, \phi_n(s,x))^T = \phi(s,x)$, and each $\phi_i(s,x)$ is bounded on $[-\tau, 0] \times \Omega_\sigma$.

In some cases, the following assumptions may be considered:

(A1) There is a positive definite diagonal matrix $G = \text{diag}(G_1, G_2, \cdots, G_n)$ such that

$|g_i(s) - g_i(t)| \leqslant G_i |s - t|, \quad \forall s, t \in \mathbb{R};$

(A2) There is a positive real number $c > 0$ such that

$$
0 \leqslant [-C_\sigma v + A_\sigma g(v) + B_\sigma g(v) + J_\sigma] \leqslant c D_\sigma E, \quad \forall v \in R^n
$$

where $D_\sigma > 0$ is a positive definite diagonal matrix, and $E = (1, 1, \cdots, 1)^T \in \mathbb{R}^n$.

**Remark 1.** The condition (A2) is so assumed that Theorem 3.1 and Theorem 3.2 can guarantee the existence of the positive stationary solution of the neural networks (2.1), for the positive stationary solution has a general adaptive property in most dynamical systems, including neural networks, ecosystem and other various real engineering systems. Very recently, Ruofeng Rao, Quanxin Zhu and Kaibo Shi in [19] generalized the conclusion of Theorem 3.1 of this paper from the Lipschitz condition to the generalized Lipschitz condition, and utilized the boundedness condition similar as (A2) to verify the existence of positive stationary solution of the Gilpin-Ayala competition model with reaction–diffusion and delayed feedback. So it is necessary that we assume the condition (A2), and give the existence theorems of positive stationary solution in this paper. In Example 4.1, the authors give the example in which the activation function $g(\cdot)$ satisfies the condition (A2). Of course, the stationary solution of neural network is not necessarily positive. So we replace (A2) with (A2∗) in Corollary 3.6–3.7 of this paper. However, even if the activation function $g(\cdot)$ is bounded, (A2) or (A2∗) can not be satisfied. On one hand, it is the profound complexity brought about by diffusion phenomenon, which makes this paper different from other related literature [2,11]. That is, diffusion has two sides. On one hand, it promotes the stability of the system; on the other hand, it brings complexity, which makes it difficult to judge the stability in many cases, even in the case of the activation function $g(\cdot)$ being bounded. This is a new viewpoint, different from those of many existing literature where only the former is emphasized.

Recently, Ruofeng Rao, Quanxin Zhu and Kaibo Shi in [19] assume that the population density of the ecosystem is limited due to the limited natural resources, and then the condition similar as (A2) is easily satisfied in [19, Theorem 3.1]. So we wish that the state variables $y(t,x)$ of neural networks are also bounded. Fortunately, in 2019, Ruofeng Rao, Shouming Zhong and Zhilin Pu utilized Laplacian semigroup theory and $L^\infty$-estimate technique to derive the boundedness theorem on the state variables of reaction–diffusion BAM neural networks [46, Theorem 3.3]. Particularly, the boundedness in [46, Theorem 3.3] is about the norm $\|\cdot\|_{L^\infty}$, which together with the continuity of the state variables implies this boundedness is similarly as that of [19,13]. Due to the limited space of this paper, we might consider this problem in the next paper (see Problem 5 for details).

Define the switching law as follows,

**Switching Law** $\mathfrak{F}$: At each switching we determine the next mode according to the following minimum law:

$$
\sigma(t) = \arg\min (y - y^\sigma)^T
$$
$$
\left[ \left( -2\lambda_{\sigma 1} D_\sigma - 2C_\sigma + A_\sigma A_\sigma^T + B_\sigma B_\sigma^T + G^2 + e^{\gamma \tau} q G^2 \right) + \Psi \right] (y - y^\sigma).
\tag{2.5}
$$

$(\mathfrak{F}_1)$ Choose the initial mode $\sigma(t) = i_0$, if $(y(t_0, x) - y^\sigma(x)) \in \Upsilon_{i_0}$.

$(\mathfrak{F}_2)$ For each $t > t_0$, if $\sigma(t^-) = i$ and $(y - y^\sigma) \in \Upsilon_i$, keep $\sigma(t) = i$. On the other hand, if $\sigma(t^-) = i$ but $(y - y^\sigma) \notin \Upsilon_i$. i.e., hitting a switching surface, choose the next mode by applying (2.5) and begin to switch.

Here, $\Psi$ is a positive definite symmetric matrix with $\lambda_{\min} \Psi > 0$, and $\Upsilon_\sigma$ is defined as follows,

$$\Upsilon_\sigma = \left\{ y \in \mathbb{R}^n | (y - y^\sigma)^T \right.$$
$$\left. \left( -2\lambda_{\sigma 1} D_\sigma - 2C_\sigma + A_\sigma A_\sigma^T + B_\sigma B_\sigma^T + G^2 + e^{\gamma\tau} q G^2 + \Psi \right)(y - y^\sigma) < 0 \right\},$$

where $\lambda_{\min}\Psi$ represents the minimum of all the eigenvalues of the symmetric matrix $\Psi > 0$.

**Definition 1.** A system is said to be globally asymptotically stable if it owns an equilibrium point which is globally asymptotically stable. Particularly, the globally asymptotical stability of an ordinary differential system is said to be invariant under the influence of diffusions if a constant equilibrium point $u^*$ of the ordinary differential system is globally asymptotically stable, and $u^*(x)$ is a globally asymptotically stable stationary solution of its corresponding reaction–diffusion system, where $u^*(x)$ is not necessarily equal to $u^*$.

**Definition 2.** ([12]) Let $\psi$ be a real $C^1$ functional defined on a Banach space $X$. If any sequence $\{u_n\}$ in $X$ with $\psi(u_n) \to a$ and $\|\psi'(u_n)\|_{X^*} \to 0$ has a convergent subsequence, and this holds for every $a \in \mathbb{R}$, one says that $\psi$ satisfies the (PS) condition.

**Definition 3.** Suppose that for each $\sigma \in \mathfrak{T}$, there exists the unique stationary solution $y^\sigma(x)$ for the switched system (2.1), and $u(t, x) = y(t, x) - y^\sigma(x)$ satisfies

$$\|u\|_{L^2(\Omega_\sigma)}^2 \leqslant M\|\phi\|_\tau^2 e^{-\gamma t}, \quad \forall t \geqslant 0,$$

where $\gamma > 0$ and $M > 1$ are constants. Then we say, the switched system (2.1) is globally exponentially stable in the meaning of switching, and the null solution of the switched delayed reaction–diffusion system (2.3) equipped with the initial value (2.4) is globally exponentially stable. Particularly, in the case of $\mathfrak{T} = \{1\}$ or $y^\sigma(x) \equiv y(x)$ for all $\sigma \in \mathfrak{T}$, we say, the switched system (2.1) is globally exponentially stable (**in the classical meaning**).

**Lemma 2.1.** ([32]) Let $\mathfrak{Z}$ be a Banach space, and $\mathfrak{R}$ is a closed convex set. If $\mathfrak{T} : \mathfrak{R} \to \mathfrak{R}$ is a compact mapping such that for any $\varphi \in \mathfrak{R}$ with $\|\varphi\| = M$, the inequality $\varphi \neq r\mathfrak{T}(\varphi)$ holds for each $r \in [0, 1]$, where $M$ is any given positive constant, then there exits at least a fixed point of $\mathfrak{T}$, say, $\varphi \in \mathfrak{R}$ with $\|\varphi\| \leqslant M$.

**Lemma 2.2.** ([12]) Let $H = H_1 \oplus H_2$ be a Banach space, and $H_1$ is a finite dimension subspace. If $\psi \in C^1(H, \mathbb{R})$, satisfying $\psi(0) = 0$, the (PS) condition. Besides, for some $\delta > 0$, the following conditions hold,

(P1) $\psi(u) \leqslant 0$ if $u \in H_1$ with $\|u\| \leqslant \delta$;
(P2) $\psi(u) \geqslant 0$ if $u \in H_2$ with $\|u\| \leqslant \delta$;
(P3) $\psi$ is bounded below, satisfying $\inf_H \psi < 0$,

then $\psi$ owns at least two non-zero critical points.

**Lemma 2.3.** ([17]). For the given matrices $E, F$, and $G$ with $F^T F \leqslant I$ and scalar $\varepsilon > 0$, the following inequality holds:

$$GFE + E^T F^T G^T \leqslant \varepsilon GG^T + \varepsilon E^T E$$

## 3. Main results

**Theorem 3.1.** Suppose that the conditions (A1) and (A2) hold, then the system (2.1) possesses a positive bounded stationary solution $y^\sigma(x)$ for $x \in \Omega_\sigma$ with $y^\sigma|_{\partial\Omega_\sigma} = 0$. In addition, there is a sequence of nonnegative constants $\beta_\sigma (\sigma = 1, 2, \cdots, N)$ with $\sum_{\sigma=1}^N \beta_\sigma = 1$ and $0 \leqslant \beta_\sigma \leqslant 1$ and positive constants $\gamma \in (0, \lambda_{\min}\Psi)$ and $q > 1$ such that

$$\sum_{\sigma=1}^N \beta_\sigma \left( -2\lambda_{\sigma 1} D_\sigma - 2C_\sigma + A_\sigma A_\sigma^T + B_\sigma B_\sigma^T \right) + G^2 + e^{\gamma\tau} q G^2 + \Psi$$
$$< 0, \tag{3.1}$$

then the null solution of the switched delayed reaction–diffusion system (2.3) equipped with the initial value (2.4) is exponentially stable with the convergence rate $\frac{\gamma}{2}$.

**Proof.** Firstly, we denote $\|u_i\| = \sqrt{\int_{\Omega_\sigma} |\nabla u_i|^2 dx}$, and $\|u\| = \sum_{i=1}^n \|u_i\|$. Besides, denote by $I$ the identity matrix. If the stationary solution of the system (2.1) exists, we may denote $y^\sigma(x)$.

Define the operator $\mathfrak{M} : \left[C\left(\overline{\Omega_\sigma}\right)\right]^n \to \left[C\left(\overline{\Omega_\sigma}\right)\right]^n$ as follows,

$$\mathfrak{M} = \begin{pmatrix} -\Delta & 0 & 0 & \cdots & 0 \\ 0 & -\Delta & 0 & \cdots & 0 \\ \vdots & \vdots & \vdots & \cdots & \vdots \\ 0 & 0 & 0 & \cdots & -\Delta \end{pmatrix}.$$

The operator $\mathfrak{M}$ has the inverse operator $\mathfrak{M}^{-1}$ as follows,

$$\mathfrak{M}^{-1} = \begin{pmatrix} (-\Delta)^{-1} & 0 & 0 & \cdots & 0 \\ 0 & (-\Delta)^{-1} & 0 & \cdots & 0 \\ \vdots & \vdots & \vdots & \cdots & \vdots \\ 0 & 0 & 0 & \cdots & (-\Delta)^{-1} \end{pmatrix},$$

where $\mathfrak{M}^{-1} : \left[C\left(\overline{\Omega_\sigma}\right)\right]^n \to \left[C\left(\overline{\Omega_\sigma}\right)\right]^n$ is a linear compact positive operator [43], and

$$\begin{cases} \mathfrak{M} y^\sigma(x) = -D_\sigma^{-1} C_\sigma y^\sigma(x) + D_\sigma^{-1} A_\sigma g(y^\sigma(x)) + D_\sigma^{-1} B_\sigma g(y^\sigma(x)) + D_\sigma^{-1} J_\sigma, & x \in \Omega_\sigma, \\ y_i^\sigma(x) = 0, x \in \partial\Omega_\sigma, i = 1, 2, \cdots, n, \end{cases}$$

It is obvious that $\left( -D_\sigma^{-1} C_\sigma y^\sigma(x) + D_\sigma^{-1} A_\sigma g(y^\sigma(x)) + D_\sigma^{-1} B_\sigma g(y^\sigma(x)) + D_\sigma^{-1} J_\sigma \right)$ is continuous for all the variables $x, y_1^\sigma, \cdots, y_n^\sigma$. Define

$$\mathfrak{R} = \left\{ \varphi(x) \in \left[ C\left(\overline{\Omega_\sigma}\right) \right]^n : \varphi(x) \geqslant 0, x \in \Omega; \varphi(x) = 0, x \in \partial\Omega; \|\varphi(x)\| < +\infty \right\},$$

then $\mathfrak{R}$ is a positive cone, which must be a closed convex subset of $\left[ C\left(\overline{\Omega_\sigma}\right) \right]^n$. Define an operator $\mathfrak{T} : \mathfrak{R} \to \mathfrak{R}$ such that

$$\mathfrak{T}\varphi = \mathfrak{M}^{-1} \left( -D_\sigma^{-1} C_\sigma \varphi + D_\sigma^{-1} A_\sigma g(\varphi) + D_\sigma^{-1} B_\sigma g(\varphi) + D_\sigma^{-1} J_\sigma \right),$$
$$\varphi \in \mathfrak{R}.$$

Because $\mathfrak{M}^{-1}$ is the linear positive compact operator [43], and $\left(-D_\sigma^{-1}C_\sigma y^\sigma(x) + D_\sigma^{-1}A_\sigma g(y^\sigma(x)) + D_\sigma^{-1}B_\sigma g(y^\sigma(x)) + D_\sigma^{-1}J_\sigma\right)$ is positive continuous, we can conclude that $\mathfrak{T} : \mathfrak{R} \to \mathfrak{R}$ is a positive compact operator.

Next, we claim that $\mathfrak{T}$ satisfies all the assumption conditions of Lemma 2.1, which implies that $\mathfrak{T}$ has at least one fixed point in $\mathfrak{R}$.

Indeed, if it is not true, there must be $\{r_n\} \subset [0,1]$ and $\{\varphi_n\} \subset \mathfrak{R}$ with

$$\varphi_n = r_n \mathfrak{T}(\varphi_n)$$
$$= r_n \mathfrak{M}^{-1}\left(-D_\sigma^{-1}C_\sigma \varphi_n + D_\sigma^{-1}A_\sigma g(\varphi_n) + D_\sigma^{-1}B_\sigma g(\varphi_n) + D_\sigma^{-1}J_\sigma\right) \tag{3.4}$$

and

$$\|\varphi_n\| = M_n \to +\infty, \quad n \to +\infty.$$

The compactness of bounded closed sets in a finite dimensional space yields that there is a subsequence of $\{r_n\}$, say, $\{r_n\}$ such that $\lim_{n\to\infty} r_n = r_0$.

Let

$$\mathfrak{L}_n = \frac{\varphi_n}{\|\varphi_n\|},$$

then it is easy to conclude from (3.4) and (A2) that if $r_n \to r_0 \in [0,1]$,

$$\mathfrak{L}_n \to \mathfrak{L}_0 \in \mathfrak{R}, \quad \|\mathfrak{L}_0\| = 1. \tag{3.5}$$

In fact, combining (H2) and the property of the operator $\mathfrak{M}^{-1}$ yields

$$\mathfrak{L}_n = r_n \mathfrak{M}^{-1}\left(\frac{-D_\sigma^{-1}C_\sigma \varphi_n + D_\sigma^{-1}A_\sigma g(\varphi_n) + D_\sigma^{-1}B_\sigma g(\varphi_n) + D_\sigma^{-1}J_\sigma}{\|\varphi_n\|}\right)$$
$$\to 0 \in R^n, \quad n \to \infty.$$

On one hand, $\mathfrak{L}_0 = 0$ implies $\|\mathfrak{L}_0\| = 0$. On the other hand, it follows by $\mathfrak{L}_n \to \mathfrak{L}_0$ and $\|\mathfrak{L}_n\| = 1$ that $\|\mathfrak{L}_0\| = 1$, which contradicts $\|\mathfrak{L}_0\| = 0$, and hence all the conditions of Lemma 2.1 are satisfied. Thereby, there exists $y^\sigma \in \mathfrak{R}$ such that $y^\sigma = \mathfrak{T}y^\sigma$ with $\|y^\sigma\| \leqslant M$, and $y^\sigma$ is a bounded positive solution of the system (2.1).

Next, we consider the following Lyapunov functional:

$$V = \int_{\Omega_\sigma} [y(t,x) - y^\sigma(x)]^T [y(t,x) - y^\sigma(x)] dx$$
$$= \int_{\Omega_\sigma} u^T(t,x) u(t,x) dx. \tag{3.8}$$

Set

$$U(t, u(t,x)) = \begin{cases} e^{\gamma t} \int_{\Omega_\sigma} u^T(t,x) u(t,x) dx, & t \geqslant 0, \\ \int_{\Omega_\sigma} u^T(t,x) u(t,x) dx, & t \in [-\tau, 0], \end{cases}$$

It is obvious that $U$ is continuous for $t \geqslant -\tau$. For $t \geqslant 0$ and $\gamma > 0$,

Now we claim that there is a positive constants $C_0 > 1$ and $K \in \mathbb{R}$ with $K > 1$ such that

$$U(t, u(t,x)) \leqslant KC_0 \|\phi\|_\tau^2, \quad \forall t \geqslant 0, \tag{3.9}$$

where $\|\phi\|_\tau^2 = \sup\limits_{s \in [-\tau, 0]} \int_{\Omega_\sigma} \phi^T(s,x) \phi(s,x) dx$.

Indeed, suppose this claim is not true, then there must be a $t \geqslant 0$ such that $U(t, u(t,x)) > KC_0 \|\phi\|_\tau^2$. Obviously, (3.9) holds for $t \in [-\tau, 0]$, and hence there must exist $t^* > 0$ such that

$$U(t^*, u(t^*,x)) = KC_0 \|\phi\|_\tau^2 \quad \text{and} \quad U(t, u(t,x)) \leqslant KC_0 \|\phi\|_\tau^2, \ \forall t \in [0, t^*],$$

and hence

$$U(t^*, u(t^*,x)) = KC_0 \|\phi\|_\tau^2 \quad \text{and} \quad U(t, u(t,x))$$
$$\leqslant KC_0 \|\phi\|_\tau^2, \ \forall t \in [-\tau, t^*]. \tag{3.10}$$

Let $q > 1$, and due to $U(0, u(0,x)) < KC_0 \|\phi\|_\tau^2 = U(t^*, u(t^*,x))$, there is $t^{**} \in [0, t^*]$ such that

$$\begin{cases} U(t^{**}, u(t^{**})) = \frac{1}{q}KC_0 \|\phi\|_\tau^2 < KC_0 \|\phi\|_\tau^2 = U(t^*, u(t^*)); \\ U(t^{**}, u(t^{**})) \leqslant U(t, u(t,x)) \leqslant U(t^*, u(t^*)) = KC_0 \|\phi\|_\tau^2, \ \forall t \in [t^{**}, t^*]. \end{cases} \tag{3.11}$$

It follows from (3.10), (3.11) and the definition of $U(t, u(t,x))$ that for $s \in [-\tau, 0]$ and $t \in [t^{**}, t^*]$,

$$\int_{\Omega_\sigma} e^{\gamma s} [u^T(t+s) u(t+s)] dx \leqslant q \int_{\Omega_\sigma} [u^T(t,x) u(t,x)] dx, \tag{3.12}$$

which yields that for any $s \in [-\tau, 0]$,

$$\int_{\Omega_\sigma} [u^T(t - \tau(t)) u(t - \tau(t))] dx \leqslant e^{\gamma \tau} q \int_{\Omega_\sigma} u^T(t,x) u(t,x) dx,$$
$$t \in [t^{**}, t^*]. \tag{3.13}$$

On the other hand, the condition (A1) yields

$$u^T A_\sigma f(u) + f^T(u) A_\sigma^T u \leqslant u^T\left(A_\sigma A_\sigma^T\right)u + u^T G^2 u,$$

and

$$u^T B_\sigma f(u(t - \tau(t), x)) + f^T(u(t - \tau(t), x)) B_\sigma^T u$$
$$\leqslant u^T\left(B_\sigma B_\sigma^T\right)u + u^T(t - \tau(t), x) G^2 u(t - \tau(t), x).$$

Now, we calculate the derivative $\frac{dV}{dt}$ alongside with the trajectories of the system (2.1) or (2.3) as follows,

$$\frac{dV}{dt} \leqslant \int_{\Omega_\sigma} u^T(t,x)\left(-2\lambda_1 D_\sigma - 2C_\sigma + A_\sigma A_\sigma^T + B_\sigma B_\sigma^T + G^2\right)u(t,x)dx$$
$$+ \int_{\Omega_\sigma} u^T(t - \tau(t), x) G^2 u(t - \tau(t), x)dx,$$

which together with (3.13) implies that

$$\frac{dV}{dt} \leqslant \int_{\Omega_\sigma} u^T(t,x)\left(-2\lambda_1 D_\sigma - 2C_\sigma + A_\sigma A_\sigma^T + B_\sigma B_\sigma^T + G^2 + e^{\gamma \tau} q G^2\right)u(t,x)dx, \ t \in [t^{**}, t^*]. \tag{3.14}$$

For any given $t \geqslant t_0$, according to the switching law $\mathfrak{F}$, when $\sigma(t^-) = i$ and $u(t,x) \in \Upsilon_i$, then keep $\sigma(t) = i$, and we can conclude that

$$\frac{dV}{dt} \leqslant \int_{\Omega_\sigma} u^T(t,x)(-\Psi) u(t,x) dx \leqslant -\lambda_{\min} \Psi \int_{\Omega_\sigma} u^T u dx$$
$$= -\lambda_{\min} \Psi V(t, u(t,x)), \ t \in [t^{**}, t^*]. \tag{3.15}$$

When $\sigma(t^-) = i$ and $u(t,x) \notin \Upsilon_i$, which means that the trajectory hits a switching surface. On the other hand, it is not difficult to deduce from (3.1) that $\bigcup\limits_{i=1}^N \Upsilon_i = R^n \setminus \{0\}$, which together with the minimum law (2.5) yields (3.15), too.

Thereby, it follows from the definition of $U(t, u(t,x))$ that

$$\frac{dU}{dt} = (\gamma - \lambda_{\min} \Psi) U(t, u(t,x)) \leqslant 0, \ t \in [t^{**}, t^*],$$

which derives that $U(t^*, u(t^*)) \leqslant U(t^{**}, u(t^{**}))$. This contradicts (3.11). So we have prove the claim (3.9), which means

$$\|u\|_{L^2(\Omega_\sigma)}^2 \leqslant KC_0 \|\phi\|_\tau^2 e^{-\gamma t}, \quad \forall t \geqslant 0,$$

which implies that the switched delayed reaction–diffusion system (2.3) equipped with the initial value (2.4) is exponentially stable with the convergence rate $\frac{\gamma}{2}$. $\quad\square$

**Remark 2.** Particularly, in the case of $\mathfrak{T} = \{1\}$ or $N = 1$, the system (2.1) is the common delayed reaction–diffusion system without any switches. Theorem 3.1 **includes** the exponential stability **in the classical sense** for the positive bounded stationary solution of the following common delayed reaction–diffusion system:

$$\begin{cases} \frac{\partial y(t,x)}{\partial t} = & D\Delta y(t,x) - Cy(t,x) + Ag(y(t,x)) + Bg(y(t-\tau(t),x)) + J, \quad (t,x) \in \mathbb{R}_+ \times \Omega, \\ y_i(t,x) = & 0, t \geqslant 0, x \in \partial\Omega, \ i = 1,2,\cdots,n. \end{cases}$$

(3.16)

**Remark 3.** For the first time, Theorem 3.1 shows the two sides of the diffusion phenomena in practical engineering (see Remark 6 for details). But [1, Theorem 1] only shows one side of the diffusion phenomena which promotes the stability of reaction–diffusion neural networks. So do those of previous related literature [2–5] and the references therein.

Next, the uniqueness of the positive stationary solution of Theorem 3.1 will be presented by adding a condition to Theorem 3.1 so that the exponential stability of the positive stationary solution is global (**in the meaning of** Definition 3).

**Theorem 3.2.** *If all the assumptions of Theorem 3.1 hold, and if, in addition, the following condition is satisfied,*

*(A3) for each $\sigma \in \mathfrak{T}$, there exists a scalar $\varepsilon > 0$ such that*

$$-C_\sigma + \frac{p_\sigma}{2}\left(\varepsilon^{-1}I + \varepsilon G^2\right) < \lambda_{\sigma 1} D_\sigma,$$

(3.17)

*where the constant $p_\sigma > 0$ satisfying $p_\sigma^2 I \geqslant (A_\sigma + B_\sigma)^T(A_\sigma + B_\sigma)$, then the system (2.1) possesses a unique positive bounded stationary solution $y^\sigma(x)$ for $x \in \Omega_\sigma$ with $y^\sigma|_{\partial\Omega_\sigma} = 0$. And the null solution of the switched delayed reaction–diffusion system (2.3) equipped with the ini-* tial value (2.4) is globally exponentially stable with the convergence rate $\frac{\gamma}{2}$. Particulary, if $\mathfrak{T} = \{1\}$ or $N = 1$, the unique stationary solution $y^\sigma(x)(\sigma = 1)$ of the **deterministic system** (2.1) is globally exponentially stable **in the classical meaning**.

**Proof.** Assume both $y(x)$ and $v(x)$ are the stationary solutions of the system (2.1). Then we claim $y(x) = v(x)$.

In fact, Lemma 2.3 yields

$$(y(x) - v(x))^T(A_\sigma + B_\sigma)(g(y(x)) - g(v(x)))$$
$$\leqslant \frac{p_\sigma}{2}(y(x) - v(x))^T\left[\varepsilon^{-1}I + \varepsilon G^2\right](y(x) - v(x)).$$

(3.18)

Since both $y(x)$ and $v(x)$ are the stationary solutions of the system (2.1), we can see it from (3.18), variational method and the Poincare inequality that

$$\lambda_{\sigma 1}\int_{\Omega_\sigma}|y(x) - v(x)|^T D_\sigma|y(x) - v(x)|dx$$
$$\leqslant \int_{\Omega_\sigma}|y(x) - v(x)|^T\left[-C_\sigma + \frac{p_\sigma}{2}\left(\varepsilon^{-1}I + \varepsilon G^2\right)\right]|y(x) - v(x)|dx.$$

Now the condition (A3) yields the claim via the proof by contradiction. And so the system (2.1) possesses a unique positive bounded stationary solution $y^\sigma(x)$ for $x \in \Omega_\sigma$ with $y^\sigma|_{\partial\Omega_\sigma} = 0$. Moreover, according to the proof of Theorem 3.1, the unique positive bounded stationary solution $y^\sigma(x)$ is globally exponentially stable, i.e., the null solution of the switched delayed reaction–diffusion system (2.3) equipped with the initial value (2.4) is globally exponentially stable with the convergence rate $\frac{\gamma}{2}$. $\quad\square$

**Discussion and comparison.** Different from previous literature related to switched system, our Theorem 3.2 presents the unique existence of the nontrivial stationary solution for the switched system so that the globally asymptotical stability of the switched system can be guaranteed.

**Remark 4.** (A2) and (A3) are the sufficient conditions, guaranteeing the global stability invariance of ordinary differential systems under the influence of diffusions in the meaning of Definition 1 (see Corollary 3.4).

To show the idea of Remark 3, we may consider the stability of the constant equilibrium point of the following delayed reaction–diffusion Cohen-Grossberg neural networks which is the partial differential equations model studied in [2]:***

$$\begin{cases} \frac{\partial u_i(t,x)}{\partial t} = & r_i\Delta u_i(t,x) - a_i(u_i(t,x))\left[b_i(u_i(t,x)) - \sum_{j=1}^{n}c_{ij}f_j(u_j(t,x)) - \sum_{j=1}^{n}d_{ij}g_j(u_j(t-\tau_j(t),x)) + I_i\right], \quad t \geqslant 0, t \neq t_k, x \in \Omega, \\ u_i(t^+,x) = & m_iu_i(t^-,x) + \sum_{j=1}^{n}n_{ij}h_j(u_j(t^--\tau_j(t),x)), \quad t = t_k, 0 \leqslant \tau_j(t) \leqslant \tau_j, \forall j, \\ u_i(t,x) = & 0, \quad t \geqslant 0, x \in \partial\Omega, i = 1,2,\cdots,n, \\ u_i(s,x) = & \phi_i(s,x), \quad -\tau \leqslant s \leqslant 0, \tau = \max_{1\leqslant i\leqslant n}\tau_j, \end{cases}$$

(3.19)

where $u_i(t_k^+,x) = u_i(t_k,x)$, all the variables, coefficients and functions are defined in [2], and are different from those of our Theorem 3.1 and Theorem 3.2. Below, we will give a stability criterion of its corresponding ordinary differential equations as follows, which will be completely similar as [2, Theorem 3.1]:

$$\begin{cases} \frac{du_i(t)}{dt} = & -a_i(u_i(t))\left[b_i(u_i(t)) - \sum_{j=1}^{n}c_{ij}f_j(u_j(t)) - \sum_{j=1}^{n}d_{ij}g_j(u_j(t-\tau_j(t))) + I_i\right], \quad t \geqslant 0, t \neq t_k, \\ u_i(t^+) = & m_iu_i(t^-) + \sum_{j=1}^{n}n_{ij}h_j(u_j(t^--\tau_j(t))), \quad t = t_k, \\ u_i(s) = & \phi_i(s), \quad -\tau \leqslant s \leqslant 0, \tau = \max_{1\leqslant j\leqslant n}\tau_j, \quad i = 1,2,\cdots,n. \end{cases}$$

(3.20)

For the convenience of readers, we may copy the assumption conditions in the document [2] as follows,

(H1) Each function $a_i(u)$ is bounded, positive and continuous, i.e., there exist two positive diagonal matrices $\underline{A} = diag(\underline{A}_1, \underline{A}_2, \cdots, \underline{A}_n)$ and $\overline{A} = diag(\overline{A}_1, \overline{A}_2, \cdots, \overline{A}_n)$ such that

$$0 < \underline{A}_i \leqslant a_i(u) \leqslant \overline{A}_i, \quad \forall u \in \mathbb{R}, \forall i.$$

(H2) Each function $b_i(u)$ is monotone increasing, i.e., there exist a positive diagonal matrix $B = diag(B_1, B_2, \cdots, B_n)$ such that

$$\frac{b_i(u) - b_i(v)}{u - v} \geqslant B_i, \quad \forall u, v(u \neq v) \in \mathbb{R}, \forall i.$$

(H3) There exist three positive diagonal matrices $F = diag(F_1, F_2, \cdots, F_n)$, $G = diag(G_1, G_2, \cdots, G_n)$ and $H = diag(H_1, H_2, \cdots, H_n)$ such that

$$0 \leqslant \frac{f_i(u) - f_i(v)}{u - v} \leqslant F_i, \quad 0 \leqslant \frac{g_i(u) - g_i(v)}{u - v} \leqslant G_i, \quad 0 \leqslant \frac{h_i(u) - h_i(v)}{u - v}$$
$$\leqslant H_i, \quad \forall u, v(u \neq v) \in \mathbb{R}, \forall i.$$

It is obvious that the system (3.19) and its corresponding ordinary differential Eqs. (3.20) own the same constant equilibrium point $u^* = (u_1^*, \cdots, u_n^*)^T \in \mathbb{R}^n$ if $u^* = (u_1^*, \cdots, u_n^*)^T$ satisfies

$$\begin{cases} b_i(u_i^*) - \sum_{j=1}^n c_{ij} f_j\left(u_j^*\right) - \sum_{j=1}^n d_{ij} g_j\left(u_j^*\right) + I_i = 0, \\ \qquad\qquad\qquad\qquad\qquad\qquad \forall i, j = 1, 2, \cdots, n. \\ (m_i - 1)u_i^* + \sum_{j=1}^n n_{ij} h_j\left(u_j^*\right) = 0, \end{cases} \quad (3.21)$$

**Theorem 3.3.** *Under assumptions (H1)–(H3), if the following conditions hold:*

$\left(\widetilde{C}1\right)$ *there exists a positive diagonal matrix $P > 0$ such that*

$$\widetilde{\Psi} = \begin{pmatrix} -2P\underline{A}B + F^2 & P\overline{A}|C| & P\overline{A}|D| \\ |C^T|\overline{A}P & -I & 0 \\ |D^T|\overline{A}P & 0 & -I \end{pmatrix} < 0,$$

*where* $R = diag(r_1, r_2, \cdots, r_n), |C| = \left(|c_{ij}|\right)_{n \times n}, \quad |D| = \left(|d_{ij}|\right)_{n \times n}, I = diag(1, 1, \cdots, 1)$;

$\left(\widetilde{C}2\right)$ $\tilde{a} = \frac{\lambda_{\min}\widetilde{\Phi}}{\lambda_{\max}P} > \frac{\lambda_{\max}G^2}{\lambda_{\min}P} = b \geqslant 0$, *where*

$\widetilde{\Phi} = 2P\underline{A}B - P\overline{A}|C||C^T|\overline{A}P - P\overline{A}|D||D^T|\overline{A}P - F^2 > 0$;

*(C3) there exists a constant $\delta$ such that $\delta > \ln(\rho e^{\lambda \tau})/\delta\tau$, where $\lambda > 0$ is the unique solution of the equation $\lambda = a - be^{\lambda \tau}$, and $\rho = \max\left\{1, \frac{2\lambda_{\max}(PMP)}{\lambda_{\min}P} + \frac{2\lambda_{\max}(HN^TPNH)}{\lambda_{\min}P} e^{\lambda \tau}\right\}, M = diag(m_1, m_2, \cdots, m_n), N = (n_{ij})_{n \times n}, H = diag(H_1, \cdots, H_n)$;*

*then we have the conclusions:*

**Conclusion** *(1) the constant equilibrium point $u^*$ of the ordinary differential system (3.20) is globally exponentially stable with convergence rate $\frac{1}{2}\left(\lambda - \frac{\ln(\rho e^{\lambda \tau})}{\delta\tau}\right)$;*

**Conclusion** *(2) If there exits a solution $u^*(x)(\neq u^*$ in common cases) of the following equations:*

$$\begin{cases} 0 = & r_i \Delta u_i(x) - a_i(u_i(x)) \left[b_i(u_i(x)) - \sum_{j=1}^n c_{ij} f_j(u_j(x)) - \sum_{j=1}^n d_{ij} g_j(u_j(x)) + I_i\right], \quad x \in \Omega \\ u_i(x) = & m_i u_i(x) + \sum_{j=1}^n n_{ij} h_j(u_j(x)), \quad \forall j, \\ u_i(x) = & 0, \quad x \in \partial\Omega, i = 1, 2, \cdots, n, \end{cases} \quad (3.22)$$

*then the globally exponential stability of the equilibrium point $u^*$ of the reaction–diffusion system (3.20) directly yields that the stationary solution $u^*(x)$ of the system (3.19) is also globally exponentially stable with the same convergence rate $\frac{1}{2}\left(\lambda - \frac{\ln(\rho e^{\lambda \tau})}{\delta\tau}\right)$, which implies that the diffusion promotes the stability, or the diffusion is not harmful to the stability.*

**Proof.** Firstly, we may prove the first conclusion of Theorem 3.1 involved in the system (3.20).

Next, for any given $i$, let $y_i(t) = u_i(t) - u_i^*$, where $u^*$ is the constant equilibrium point of the system (3.20). Then the system (3.20) can be transformed into

$$\begin{cases} \frac{dy_i(t)}{dt} = & -\tilde{a}_i(y_i(t))\left[\tilde{b}_i(y_i(t)) - \sum_{j=1}^n c_{ij}\tilde{f}_j(y_j(t)) - \sum_{j=1}^n d_{ij}\tilde{g}_j(y_j(t - \tau_j(t)))\right], t \geqslant 0, t \neq t_k, \\ y_i(t^+) = & m_i y_i(t^-) + \sum_{j=1}^n n_{ij}\tilde{h}_j(y_j(t - \tau_j(t))), t = t_k, \\ y_i(s) = & \phi_i(s) - u_i^*(s), -\tau \leqslant s \leqslant 0, \tau = \max_{1 \leqslant j \leqslant n}\tau_j, i = 1, 2, \cdots, n. \end{cases} \quad (3.23)$$

where $\tilde{a}_i(y_i(t)) = a_i(y_i(t) + u_i^*), \tilde{b}_i(y_i(t)) = b_i(y_i(t) + u_i^*) - b_i(u_i^*), \tilde{f}_j(y_j(t)) = f_j\left(y_j(t) + u_j^*\right) - f_j\left(u_j^*\right), \tilde{g}_j(y_j(t)) = g_j\left(y_j(t) + u_j^*\right) - g_j\left(u_j^*\right), \tilde{h}_j(y_j(t)) = h_j\left(y_j(t) + u_j^*\right) - h_j\left(u_j^*\right)$ for all $i, j = 1, 2, \cdots, n$.

Similarly as the proof of [2, Theorem 3.1], we may set up the Lyapunov function as follows,

$$\widetilde{V}(t) = Y^T(t)PY(t) = |Y^T(t)|P|Y(t)|,$$

where $Y(t) = (y_1(t), \cdots, y_n(t))^T, P = diag(p_1, p_2, \cdots, p_n)$.

For the case of $t \neq t_k$, we compute the Dini derivative of $\widetilde{V}(t)$ alongside with the trajectories of (3.23),

$$\begin{aligned} D^+\widetilde{V}(t) = & -2Y^T(t)P\widetilde{A}(Y(t))\widetilde{B}(Y(t)) + 2Y^T(t)P\widetilde{A}(Y(t))C\widetilde{F}(Y(t)) \\ & + 2Y^T(t)P\widetilde{A}(Y(t))D\widetilde{G}(Y(t - \tau(t))) \\ \leqslant & -|Y^T(t)|\widetilde{\Phi}|Y(t)| + Y^T(t - \tau(t))G^2 Y(t - \tau(t)) \\ \leqslant & -\tilde{a}\widetilde{V}(t) + b\left[\widetilde{V}(t)\right]_\tau, \end{aligned}$$

where $\widetilde{A}(Y(t)) = diag(\tilde{a}_1(y_1(t)), \cdots, \tilde{a}_n(y_n(t)))$,
$\widetilde{B}(Y(t)) = \left(\tilde{b}_1(y_1(t)), \cdots, \tilde{b}_n(y_n(t))\right)^T$,
$\widetilde{F}(Y(t)) = \left(\tilde{f}_1(y_1(t)), \cdots, \tilde{f}_n(y_n(t))\right)^T$,
$\widetilde{G}(Y(t)) = (\tilde{g}_1(y_1(t)), \cdots, \tilde{g}_n(y_n(t)))^T$,
$\widetilde{G}(Y(t - \tau(t))) = (\tilde{g}_1(y_1(t - \tau_1(t))), \cdots, \tilde{g}_n(y_n(t - \tau_n(t))))^T$.
When $t = t_k$, using the similar methods in the proof of [2, Theorem 3.1] results in that

$$\widetilde{V}(t_k) = Y^T(t_k)PY(t_k)$$
$$\leqslant 2\frac{\lambda_{\max}(PMP)}{\lambda_{\min}P}\widetilde{V}(t_k^-) + 2\frac{\lambda_{\max}\left(HN^TPNH\right)}{\lambda_{\min}P}\left[\widetilde{V}(t_k^-)\right]_\tau.$$

Now it follows from $\left(\widetilde{C}1\right), \left(\widetilde{C}2\right)$, (C3) and [2, Lemma 2.2] that

$$\widetilde{V}(t) \leqslant \rho\left[\widetilde{V}(0)\right]_\tau e^{-\left(\lambda - \frac{\ln\left(\rho e^{\lambda\tau}\right)}{\delta\tau}\right)t}, \quad t \geqslant 0,$$

or

$$\sqrt{(u(t) - u^*)^T(u(t) - u^*)} = \sqrt{Y^T(t)Y(t)} \leqslant \sqrt{\frac{\rho\lambda_{\max}P}{\lambda_{\min}P}}\sqrt{\left[Y^T(0)Y(0)\right]_\tau}$$

$$e^{-\frac{1}{2}\left(\lambda - \frac{\ln\left(\rho e^{\lambda\tau}\right)}{\delta\tau}\right)t}, \quad t \geqslant 0,$$

which has proved that the equilibrium point $u^*$ of system (3.20) is globally exponentially stable with convergence rate $\frac{1}{2}\left(\lambda - \frac{\ln(\rho e^{\lambda\tau})}{\delta\tau}\right)$.

Finally, we shall prove the second conclusion of Theorem 3.3.

Indeed, let $Y(t,x) = u(t,x) - u^*(x)$, where the stationary solution $u^*(x) = (u_1^*(x), \cdots, u_n^*(x))^T$, then the system (3.19) can be transformed into

$$\mathcal{V}(t_k) = \int_\Omega Y^T(t_k,x)PY(t_k,x)dx$$

$$\leqslant 2\frac{\lambda_{\max}(PMP)}{\lambda_{\min}P}\mathcal{V}(t_k^-) + 2\frac{\lambda_{\max}\left(HN^TPNH\right)}{\lambda_{\min}P}\left[\mathcal{V}(t_k^-)\right]_\tau,$$

and

$$\mathcal{V}(t) \leqslant \rho[\mathcal{V}(0)]_\tau e^{-\left(\lambda - \frac{\ln\left(\rho e^{\lambda\tau}\right)}{\delta\tau}\right)t}, \quad t \geqslant 0,$$

or

$$\sqrt{\int_\Omega(u(t,x) - u^*)^T(u(t,x) - u^*)dx}$$

$$\leqslant \sqrt{\frac{\rho\lambda_{\max}P}{\lambda_{\min}P}}\sqrt{\left[\int_\Omega Y^T(0,x)Y(0,x)dx\right]_\tau}\, e^{-\frac{1}{2}\left(\lambda - \frac{\ln\left(\rho e^{\lambda\tau}\right)}{\delta\tau}\right)t}, \quad t \geqslant 0,$$

which has proved that the equilibrium point $u^*$ of system (3.20) is globally exponentially stable with convergence rate $\frac{1}{2}\left(\lambda - \frac{\ln(\rho e^{\lambda\tau})}{\delta\tau}\right)$.

$$\begin{cases} \frac{\partial y_i(t,x)}{\partial t} = & r_i\Delta y_i(t,x) - \tilde{a}_i(y_i(t,x))\left[\tilde{b}_i(y_i(t,x)) - \sum_{j=1}^n c_{ij}\tilde{f}_j(y_j(t,x)) - \sum_{j=1}^n d_{ij}\tilde{g}_j(y_j(t - \tau_j(t),x))\right], \quad t \geqslant 0, t \neq t_k, \\ y_i(t^+,x) = & m_i y_i(t^-,x) + \sum_{j=1}^n n_{ij}\tilde{h}_j(y_j(t^- - \tau_j(t),x)), \quad t = t_k, 0 \leqslant \tau_j(t) \leqslant \tau_j, \forall j, \\ y_i(t,x) = & 0, \quad t \geqslant 0, x \in \partial\Omega, i = 1, 2, \cdots, n, \\ y_i(s,x) = & \phi_i(s,x) - u_i^*(x), \quad -\tau \leqslant s \leqslant 0, \tau = \max_{1 \leqslant j \leqslant n}\tau_j, \end{cases} \tag{3.24}$$

where $y_i, \tilde{a}_i, \tilde{b}_i, \tilde{f}_j, \tilde{g}_j$ and $\tilde{h}_i$ all are defined as those of [2].

$$\int_\Omega |\nabla u_i(t,x)|^2 dx \geqslant \lambda_1 \int_\Omega u_i^2(t,x)dx$$

On the other hand, the Poincare inequality and the Dirichlet zero boundary value yields

$$\int_\Omega Y^T(t,x)PR\Delta Y(t,x)dx = -\int_\Omega \sum_{i=1}^n p_i r_i \sum_{j=1}^m \left(\frac{\partial y_i}{\partial x_j}\right)^2 dx$$

$$\leqslant -\lambda_1 \int_\Omega \sum_{i=1}^n p_i r_i y_i^2(t,x)dx \leqslant 0, \tag{3.25}$$

where $\lambda_1$ is the smallest positive eigenvalue of the following eigenvalue problem:

$$\begin{cases} -\Delta\varphi(x) = \lambda\varphi(x), \quad x \in \Omega \subset \mathbb{R}^m, \\ \varphi(x) = 0, \quad x \in \partial\Omega. \end{cases}$$

Constructing the Lyapunov functional as follows,

$$\mathcal{V}(t) = \int_\Omega Y^T(t,x)PY(t,x)dx = \int_\Omega |Y^T(t,x)|P|Y(t,x)|dx$$

For the case of $t \neq t_k$, we compute the Dini derivative of $\mathcal{V}(t)$ alongside with the trajectories of (3.24),

In fact, due to (3.25), the globally exponential stability of the equilibrium point $u^*$ of system (3.20) directly yields that the equilibrium point $u^*(x)$ of system (3.19) is also globally exponentially stable with the same convergence rate $\frac{1}{2}\left(\lambda - \frac{\ln(\rho e^{\lambda\tau})}{\delta\tau}\right)$, which implies that the diffusion promotes the stability. The proof is completed. $\square$

**Remark 5.** Obviously, [2, Theorem 3.1] is the direct corollary of the conclusion (1) of Theorem 3.3 due to the Poincare inequality. In fact, the condition (C1) of [2, Theorem 3.1] is as follows,

$$\Psi = \widetilde{\Psi} + \begin{pmatrix} -2lPR & 0 & 0 \\ 0 & 0 & 0 \\ 0 & 0 & 0 \end{pmatrix} < 0 \tag{3.26}$$

and the condition (C2) of [2, Theorem 3.1] is as follows,

$$a = \frac{\lambda_{\min}\Phi}{\lambda_{\max}P} \geqslant \tilde{a} = \frac{\lambda_{\min}\widetilde{\Phi}}{\lambda_{\max}P} > \frac{\lambda_{\max}G^2}{\lambda_{\min}P} = b \geqslant slant 0, \tag{3.27}$$

where

$$\Phi = 2lPR + \widetilde{\Phi} > 0. \tag{3.28}$$

$$D^+\mathcal{V}(t) = \int_\Omega\left[-2Y^T(t,x)P\widetilde{A}(Y(t,x))\widetilde{B}(Y(t,x)) + 2Y^T(t,x)P\widetilde{A}(Y(t,x))C\widetilde{F}(Y(t,x)) + 2Y^T(t,x)P\widetilde{A}(Y(t,x))D\widetilde{G}(Y(t - \tau(t),x))\right]dx$$

$$\leqslant \int_\Omega\left[-|Y^T(t,x)|\widetilde{\Phi}|Y(t,x)| + Y^T(t - \tau(t),x)G^2Y(t - \tau(t),x)\right]dx \leqslant \int_\Omega\left[-\tilde{a}\mathcal{V}(t,x) + b[\mathcal{V}(t,x)]_\tau\right]dx,$$

where $\widetilde{A}, \widetilde{B}, \widetilde{F}, \widetilde{G}$ all are defined as those of [2]. Completely similar as the proof of the first conclusion of Theorem 3.3, we can also obtain

**Remark 6.** Conclusion (1) of Theorem 3.3 involved in the ordinary differential system (3.20) is completely similar as [2, Theorem 3.1] of the reaction–diffusion system (3.19). Due to Poincare inequality, [2, Theorem 3.1] becomes actually a corollary of the Conclusion (1) of Theorem 3.3, which implies that diffusions only promote the stability of the reaction–diffusion system. Actually, [2, Theorem 3.1] does not illuminate any negative effects on the stability of the reaction–diffusion system, compared with the conclusion (1) of Theorem 3.3. However, Conclusion (2) of Theorem 3.3 can show the two sides of the influence of diffusion on judging the stability of the diffusion system (3.19), for the existence of the solution $u^*(x)$ of the Eqs. (3.22) may place more restrictions on the system, which may be similar as the condition (A2) of Theorem 3.1.

**Remark 7.** Particularly let $a_i(u_i) \equiv 1$ and $b_i(u_i) = b_i u_i$ with $b_i \in \mathbb{R}$ in the delayed reaction–diffusion Cohen-Grossberg neural networks (3.19), then the Cohen-Grossberg neural networks (3.19) is reduced to the following cellular neural networks

$$
\begin{cases}
\frac{\partial u_i(t,x)}{\partial t} = r_i \Delta u_i(t,x) - b_i u_i(t,x) + \sum_{j=1}^{n} c_{ij} f_j(u_j(t,x)) + \sum_{j=1}^{n} d_{ij} g_j(u_j(t - \tau_j(t),x)) - I_i, & t \geqslant 0, t \neq t_k, \\
u_i(t^+,x) = m_i u_i(t^-,x) + \sum_{j=1}^{n} n_{ij} h_j(u_j(t^- - \tau_j(t),x)), & t = t_k, 0 \leqslant \tau_j(t) \leqslant \tau_j, \forall j, \\
u_i(t,x) = 0, & t \geqslant 0, x \in \partial\Omega, i = 1, 2, \cdots, n, \\
u_i(s,x) = \phi_i(s,x), & -\tau \leqslant s \leqslant 0, \tau = \max_{1 \leqslant j \leqslant n} \tau_j.
\end{cases}
$$

$$
\begin{cases}
\frac{\partial u_i(t,x)}{\partial t} = r_i \Delta u_i(t,x) - b_i u_i(t,x) + \sum_{j=1}^{n} c_{ij} g_j(u_j(t,x)) + \sum_{j=1}^{n} d_{ij} g_j(u_j(t - \tau_j(t),x)) - I_i, & t \geqslant 0, x \in \Omega, \\
u_i(t,x) = 0, & t \geqslant 0, x \in \partial\Omega, i = 1, 2, \cdots, n, \\
u_i(s,x) = \phi_i(s,x), & -\tau \leqslant s \leqslant 0, \tau = \max_{1 \leqslant j \leqslant n} \tau_j.
\end{cases}
$$

So the conclusions of Theorem 3.3 include the case of cellular neural networks. But if there is not impulse control in Theorem 3.3, **the unique existence of the stationary solution of reaction–diffusion system should be given so that the global stability can be guaranteed for the reaction–diffusion system.**

The following corollary can be derived by our Theorem 3.2 and Theorem 3.3.

**Corollary 3.4. (Global Stability Invariance).** *Suppose the conditions* (H1)–(H3), $(\overline{C}1)$ *and* $(\overline{C}2)$ *are satisfied, and* $a_i(u_i) \equiv 1, b_i(u_i) = b_i u_i, f_i(u_i) = g_i(u_i)$ *as said in* Remark 7. *Besides,* $H_i \equiv 0, m_i \equiv 1$, $n_{ij} \equiv 0, \underline{A} = \overline{A} = I$.

then we have the conclusions:

**Conclusion** (1) the constant equilibrium point $u^*$ of the following system is globally exponentially stable:

$$
\begin{cases}
\frac{du_i(t)}{dt} = -b_i u_i(t) + \sum_{j=1}^{n} c_{ij} g_j(u_j(t)) + \sum_{j=1}^{n} d_{ij} g_j(u_j(t - \tau_j(t))) - I_i, & t \geqslant 0, \\
u_i(s) = \phi_i(s), & -\tau \leqslant s \leqslant 0, \tau = \max_{1 \leqslant j \leqslant n} \tau_j, \ i = 1, 2, \cdots, n.
\end{cases}
\tag{3.29}
$$

**Conclusion** (2) If there exits a solution $u^*(x) (\neq u^*$ in common cases) of the following equations:

$$
\begin{cases}
0 = r_i \Delta u_i(x) - b_i u_i(x) + \sum_{j=1}^{n} c_{ij} g_j(u_j(x)) + \sum_{j=1}^{n} d_{ij} g_j(u_j(x)) - I_i, & x \in \Omega \\
u_i(x) = 0, & x \in \partial\Omega, i = 1, 2, \cdots, n,
\end{cases}
\tag{3.30}
$$

and if, in addition, there exists a scalar $\varepsilon > 0$ such that

$$
-B + \frac{p}{2} \left( \varepsilon^{-1} I + \varepsilon G^2 \right) < \lambda_1 R,
\tag{3.31}
$$

where the constant $p > 0$ satisfying $p^2 I \geqslant (C + D)^T (C + D)$, and $R = diag(r_1, r_2, \cdots, r_n), B = diag(b_1, b_2, \cdots, b_n), C = (c_{ij}), D = (d_{ij})$, then the globally exponential stability of the equilibrium point $u^*$ of the system (3.29) directly yields that the equilibrium point $u^*(x)$ of the following reaction–diffusion system is also globally exponentially stable (**in the classical meaning**):

**Remark 8.** Corollary 3.4 illuminates the invariance of global stability in the meaning of Definition 1. Without the impulse control, the unique existence of the stationary solution must be considered for the reaction–diffusion system. In the conclusion (2) of Corollary 3.4, the condition (3.31) is the condition (A3). Besides, the existence of the stationary solution $u^*(x)$ of the Eqs. (3.30) might be guaranteed by some condition similarly as (A2). In common cases, $u^*(x)$ is not necessarily equal to $u^*$. In fact, if the constant vector $u^* \neq 0, u^*$ is not any stationary solutions of reaction–diffusion system under Dirichlet zero boundary value. Below, an example will be designed to show it (see Statement 1).

To illuminates the effectiveness of Corollary 3.4, the Global Stability Invariance, we may present the following simple example in the case of $n = 1$.

**Example 3.5.** In Corollary 3.4, set $n = 1$, and $B = 2, C = 0.01, D = 0.01, I_1 = 0.1, g(u) = u$ and then the Lipschitz constant of $g$ is $G = 1$. Let $u^* = \frac{0.1}{1.98}$ is the unique solution of the following equation:

$$
0 = -2u + 0.01u + 0.01u + 0.1,
$$

where $u$ represents an unknown number of the above algebraic equation.

Moreover, it is easy to verify that the conclusion (1) of Corollary 3.4 holds due to the fact that the related conditions of Corollary 3.4 are satisfied. That is, the constant equilibrium point $u^* = \frac{0.1}{1.98}$ of the system (3.29) is globally exponentially stable.

Besides, set $R = 0.1, p = 0.02, \varepsilon = 1$, then $-B + \frac{p}{2}\left(\varepsilon^{-1}I + \varepsilon G^2\right) < 0 < \lambda_1 R$, where $I = 1$. According to Corollary 3.4, if $u^*(x)$ is a solution of the following equation

$$0 = 0.1\Delta u - 2u + 0.01u + 0.01u + 0.1, x \in \Omega; \quad u = 0, \quad x \in \partial\Omega,$$

whose solution is corresponding to the critical point of the following functional

$$\chi(u) = \frac{1}{2}\int_\Omega |\nabla u|^2 dx + \frac{1}{2}\int_\Omega 19.8u^2 dx - \int_\Omega u dx,$$

where $u = u(x)$, independent of $t$, then the unique stationary solution $u^*(x)$ of the system (3.30) is globally exponentially stable **in the classical meaning**. As the selected special example, we are willing to prove the existence of the stationary solution $u^*(x)$. In fact, $\chi \in C^1\left(H_0^1(\Omega), \mathbb{R}\right)$, and $\chi$ is coercive, for

$$\chi(u) = \frac{1}{2}\int_\Omega |\nabla u|^2 dx + \frac{1}{2}\int_\Omega 19.8u^2 dx - \int_\Omega u dx \geqslant \frac{1}{2}\|u\|^2 - c^*\|u\|$$
$$\to +\infty, \quad \|u\| \to \infty,$$

where $\|\cdot\|$ is the norm of $H_0^1(\Omega)$ with $\|u\|^2 = \int_\Omega |\nabla u|^2 dx$, and $c^* > 0$ is a constant. It is easy to prove that $\chi$ satisfy the Palais–Smale condition, and $\chi$ is bounded below (see the methods used in proof of Statement 2 below). And hence, $\chi$ can attain its global minimum, say, $\chi(u^*(x))$ at the point $u^*(x)$, on $H_0^1(\Omega)$ (see, e.g. the proof of [18, Lemma 2.5]). Due to the condition (3.31), the unique stationary solution $u^*(x)$ of the system (3.30) is globally exponentially stable **in the classical meaning**. Particularly, $u^*(x) \neq \frac{0.1}{1.98}, x \in \Omega$, for $u^* = \frac{0.1}{1.98}$ is a non-zero constant (see Statement 1 for details).

**Statement 1.** Let $u^*$ be a non-zero constant equilibrium point of an ordinary differential system. Then $u^*$ is not any stationary solutions of its corresponding reaction–diffusion system under Dirichlet zero boundary value, or it must lead to a contradiction.

$$\begin{cases} \frac{\partial u(t,x)}{\partial t} = & 0.001\Delta u(t,x) - 1.8u(t,x) + 0.01u(t,x) + 0.005u(t-\tau(t),x) + 1, \quad t \geqslant 0, x \in \Omega = (0,1), \\ u(t,0) = & u(t,1) = 0, \\ u(s,x) = & \xi(s,x) \text{ is bounded in } [-\tau, 0] \times (0,1), \end{cases}$$

**Proof.** Consider the following cellular neural networks in the case of $n = 1$:

$$\begin{cases} \frac{dx(t)}{dt} = & -Cx(t) + Af(x(t)) + Bf(x(t-\tau(t))) + J, \quad t \geqslant 0, \\ x(s) = & \xi(s) \text{ is bounded in } [-\tau, 0]. \end{cases} \tag{3.32}$$

and its corresponding reaction–diffusion system:

$$\begin{cases} \frac{\partial u(t,x)}{\partial t} = & D\Delta u(t,x) - Cu(t,x) + Af(u(t,x)) + Bf(u(t-\tau(t),x)) + J, \\ & (t,x) \in \mathbb{R}_+ \times \Omega, \ \Omega = (0,1) \subset \mathbb{R}, \\ u(t,x) = & 0, t \geqslant 0, x \in \partial\Omega, \\ u(s,x) = & \xi(s,x) \text{ is bounded in } [-\tau, 0] \times (0,1). \end{cases} \tag{3.33}$$

where $D = 0.001, C = 1.8, A = 0.2, B = 0.1, J = 1.09$, and $f(s) = 0.05(s-6)$ for $s \in \mathbb{R}^1$.

Direct computation derives that $x = \frac{1000}{1785} = \frac{200}{357}$ is the constant equilibrium point of the system (3.32). Obviously, $f$ is Lipschitz continuous function with Lipschitz constant $\bar{l} = 0.05$ . In [31, Theorem 1], let $p = 1$ and $N = 1$, then

$$-\mu_1\left(-NCN^{-1}\right) - \bar{l}\|NA\|_1\|N^{-1}\|_1$$

and hence [31, Theorem 1] results in that there is the unique constant equilibrium point $x^* = \frac{1000}{1785} = \frac{200}{357}$ of the ordinary differential system (3.32) is globally exponentially stable.

Next, direct computation illuminates that the stationary solutions of the system (3.33) satisfies the following equation:

$$\begin{cases} \frac{d^2u(x)}{dx^2} = 1785u(x) - 1000, \quad x \in \Omega = (0,1), \\ u(0) = u(1) = 0. \end{cases} \tag{3.34}$$

Now it is easy to verify that

$$u_*(x) = \frac{200\left(e^{-\sqrt{1785}} - 1\right)}{357\left(e^{\sqrt{1785}} - e^{-\sqrt{1785}}\right)}e^{\sqrt{1785}x} - \frac{200\left(e^{\sqrt{1785}} - 1\right)}{357\left(e^{\sqrt{1785}} - e^{-\sqrt{1785}}\right)}$$
$$e^{-\sqrt{1785}x} + \frac{200}{357}, \quad \forall x \in [0,1].$$
$$\tag{3.35}$$

is a solution of the Eq. (3.34), and it is also a nontrivial stationary solution of the reaction–diffusion (3.33).

Below, by the proof by contradiction, we shall prove that $u^*$ is not the equilibrium point of the reaction–diffusion (3.33), where $u^* = u^*(x) \equiv \frac{1000}{1785} = \frac{200}{357}$ for all $x \in (0,1)$ with $u^*(0) = u^*(1) = 0$.

Firstly, by direct computation, the system (3.33) becomes the following system:

which is equivalent to the following system via the transformation $y(t,x) = u(t,x) - u^*$:

$$\begin{cases} \frac{\partial y(t,x)}{\partial t} = & 0.001\Delta y(t,x) - 1.8y(t,x) + 0.01y(t,x) + 0.005y(t-\tau(t),x), \quad t \geqslant 0, x \in \Omega = (0,1), \\ y(t,0) = & y(t,1) = 0, \\ y(s,x) = & \eta(s,x) \text{ is bounded in } [-\tau, 0] \times (0,1), \end{cases} \tag{3.36}$$

Consider the Lyapunov functional: $V(t,y(t,x)) = \int_\Omega y^2(t,x)dx$, then the derivative $\frac{dV}{dt}$ alongside with the trajectories of the system (3.36) yields

$$\frac{dV(t,y(t,x))}{dt} = \int_\Omega 2y(t,x)(0.001\Delta y(t,x) - 1.8y(t,x) + 0.01y(t,x)$$
$$+ 0.005y(t-\tau(t),x))dx \leqslant -(0.002\pi^2 + 3.575)\int_\Omega y^2(t,x)dx$$
$$+ 0.005\int_\Omega y^2(t-\tau(t),x)dx = -aV(t,y(t,x)) + bV(t,y(t-\tau(t),x)),$$

where $a = 0.002\pi^2 + 3.575$, $b = 0.005$, satisfying $a > b > 0$. By employing [44, Lemma 3] or the methods in the proof of [13, Theorem 3], we can derive that the zero solution of the system (3.43) is globally exponentially stable with the convergence rate $\frac{\lambda}{2}$, where $\lambda > 0$ is the unique solution of the equation $\lambda = a - be^{\varsigma\tau}$. That is, the constant equilibrium point $u^*$ of the reaction diffusion system (3.33) is globally exponentially stable, and so $u^*$ is the unique equilibrium point of the reaction–diffusion system. However, $u_*(x)$ defined as (3.33) is its another equilibrium point. This contradiction shows that $u^*$ can not be the stationary solution of the reaction–diffusion system. □

**Remark 9.** In [2, Theorem 3.1], $u^*$ is a non-zero constant vector in common cases. But [2, Theorem 3.1] told us that $u^*$ is a stable equilibrium point of a reaction–diffusion system under Dirichlet zero boundary value. Now Statement 1 illuminates that it must lead to a contradiction. On the other hand, the conclusion of Theorem 3.3 points out that $u^*(x)$ a stable equilibrium point of the reaction–diffusion system under Dirichlet zero boundary value, but $u^*(x)$ is not equal to the non-zero constant vector $u^*$. This viewpoint is different from those of many previous literature (see, e.g. [2,11] and the related references therein).

**Remark 10.** Of course, some suitable non-zero constant vectors can be the equilibrium points or stationary solutions of delayed reaction–diffusion systems under Neumann boundary value (see, e.g. [9,10,45]) though they can not be the equilibrium points or stationary solutions of delayed reaction–diffusion systems under Dirichlet zero boundary value.

In the proof of Statement 1, $u^*$ is the unique equilibrium point of ordinary differential system with Lipschitz assumption on active function $f$. Now we want to know whether the number of equilibrium points changes under the influence of inevitable diffusions.

Consider the following cellular neural networks in the case of $n = 1$,

$$\frac{dx(t)}{dt} = -Cx(t) + Af(x(t)) + Bf(x(t - \tau(t))) + J, \quad \text{and } x$$
$$\in \mathbb{R}^1, \tag{3.37}$$

and its corresponding reaction–diffusion cellular neural networks

$$\begin{cases} \frac{\partial u(t,x)}{\partial t} = & D\Delta u(t,x) - Cu(t,x) + Af(u(t,x)) + Bf(u(t - \tau(t),x)) + J, \\ u(t,x) = & 0, \quad x \in \partial\Omega, \end{cases} \tag{3.38}$$

where $\Omega$ is an open bounded domain in $\mathbb{R}^3$ with smooth boundary $\partial\Omega$, $D \in \mathbb{R}^1$ is the diffusion coefficient with $D > 0$, and $C$, $A$ both are positive real numbers, $J = 0$, $B = 0$, the function $f$ is defined as follow,

$$f(u) = \begin{cases} \frac{3D}{A}\mu_1 u^{\frac{1}{3}} + \frac{2D}{A}\mu_1, & u \leqslant -1; \\ \frac{D}{A}\mu_1 u, & u \in [-1,1]; \\ \frac{3D}{A}\mu_1 u^{\frac{1}{3}} - \frac{2D}{A}\mu_1, & u \geqslant 1. \end{cases} \tag{3.39}$$

Here, we denote by $\mu_i$ the $i$th positive eigenvalue of the following eigenvalue problem:

$$\begin{cases} -\Delta u(x) + \frac{C}{D}u(x) = & \mu u(x), x \in \Omega, \\ u(x) = & 0, \quad u \in \partial\Omega, \end{cases} \tag{3.40}$$

then $\mu_1 = \frac{C}{D} + \lambda_1$, and $\mu_2 > \mu_1$ [33].

**Statement 2.** If zero solution is the global stable unique equilibrium point of ordinary differential system (3.37), then its corresponding reaction–diffusion system (3.38) owns zero solution and other stationary solutions which are at least two non-zero functions or infinitely many positive functions and negative functions.

**Proof.** Firstly, it is easy to see from $f(0) = 0$ and $J = 0$ that zero solution is also an equilibrium point of the system (3.38).

Besides, we claim that the system (3.38) owns other stationary solutions which are at least two non-zero functions or infinitely many positive functions and negative functions.

In fact, we know from (3.39) that $|\frac{df(u)}{du}| \leqslant \frac{D}{A}\mu_1$ for all $u \in \mathbb{R}^1$, because

$$\frac{df(u)}{du} = \begin{cases} \frac{D}{A}\mu_1 u^{-\frac{2}{3}}, & u \leqslant -1; \\ \frac{D}{A}\mu_1, & u \in [-1,1]; \\ \frac{D}{A}\mu_1 u^{-\frac{2}{3}}, & u \geqslant 1. \end{cases} \tag{3.41}$$

And hence, $f$ is Lipschitz continuous as follow,

$$|f(u) - f(v)| \leqslant \frac{D}{A}\mu_1 |u - v|, \quad \forall u, v \in \mathbb{R}^1. \tag{3.42}$$

And the definition of $f$ yields

$$F(u) = \int_0^u f(s)ds = \begin{cases} \frac{9}{4}\frac{D}{A}\mu_1 u^{\frac{4}{3}} + \frac{2D}{A}\mu_1 u + \frac{1}{4}\frac{D}{A}\mu_1, & u \leqslant -1; \\ \frac{1}{2}\frac{D}{A}\mu_1 u^2, & u \in [-1,1]; \\ \frac{9}{4}\frac{D}{A}\mu_1 u^{\frac{4}{3}} - \frac{2D}{A}\mu_1 u + \frac{1}{4}\frac{D}{A}\mu_1, & u \geqslant 1. \end{cases} \tag{3.43}$$

Besides, if $u(x)$ is the stationary solution of the system (3.38), $u(x)$ is a solution of the following equation:

$$\begin{cases} 0 = & D\Delta u(x) - Cu(x) + Af(u(x)), \ x \in \Omega, \\ u(x) = & 0, \quad x \in \partial\Omega, \end{cases}$$

whose solution is corresponding to the critical point of the following function:

$$\mathfrak{I}(u) = \frac{1}{2}\int_\Omega |\nabla u|^2 dx + \frac{C}{2D}\int_\Omega u^2 dx - \frac{A}{D}\int_\Omega F(u)dx,$$

where $F(u)$ is given by (3.43).

Denote $H = \left\{ u \in W_0^{1,2}(\Omega), \int_\Omega |\nabla u|^2 dx + \frac{C}{D}\int_\Omega u^2 dx < \infty \right\}$, in which the inner product is presented as follows,

$$\langle u, v \rangle = \int_\Omega \nabla u \nabla v dx + \frac{C}{D}\int_\Omega u v dx,$$

and its induced norm is denote by $\| \cdot \|_H$. Obviously, $\mathfrak{I} \in C^1(H, \mathbb{R}^1)$.

Besides, we claim that the functional $\mathfrak{I}$ is bounded below. In fact, if it is not true, there must exist a sequence $\{u_n\}$ in $H$ such that $\mathfrak{I}(u_n) \to -\infty$ as $\|u_n\|_H \to \infty$. So there must exist a scalar $\mathfrak{P} > 0$ such that $\mathfrak{I}(u_n) \leqslant \mathfrak{P}$ for all $n$. Set $\mathfrak{W}_n = u_n / \|u_n\|_H$, then $\|\mathfrak{W}_n\|_H = 1$ for all $n$. Moreover, Sobolev embedding theorem tells us that there exists $\mathfrak{W} \in H$ such that $\mathfrak{W}_n \rightharpoonup \mathfrak{W}$ in $H$, $\mathfrak{W}_n \to \mathfrak{W}$ in $L^q(\Omega)$ with $2 \leqslant q < 2^*$, and $\mathfrak{W}_n(x) \to \mathfrak{W}(x)$, a.e. $x \in \Omega$, where $2^*$ is the critical sobolev exponent. **Here and below, a subsequence of $\{u_n\}$ is still denoted by $\{u_n\}$** for convenience.

By employing the similar methods used in the proof of [18, Lemma 2.1], we claim that $\|\mathfrak{W}\|_H^2 = \mu_1 \|\mathfrak{W}\|_{L^2(\Omega)}^2$.

In fact, on one hand, it follows by (3.43) and the methods used in the proof of [18, Lemma 2.1] that there is such a constant $M > 1$ big enough that

$$\frac{\mathfrak{P}}{\|u_n\|_H^2} \geqslant \frac{1}{2}\left[ \|\mathfrak{W}_n\|_H^2 - \frac{1}{2}\mu_1 \|\mathfrak{W}_n\|_{L^2(\Omega)}^2 \right] + 0 + \frac{\int_{1 \leqslant |u_n| \leqslant M} k_M dx}{\|u_n\|_H^2} + 0.$$

On the other hand,

$|\frac{\int_{1\leqslant|u_n|\leqslant M}k_M dx}{\|u_n\|_H^2}|\leqslant\frac{|k_M|\cdot mes\Omega}{\|u_n\|_H^2}\to 0,\quad n\to\infty.$

Hence,

$\limsup_{n\to\infty}\|\mathfrak{W}_n\|_H^2-\mu_1\|\mathfrak{W}\|_{L^2(\Omega)}^2\leqslant 0\Rightarrow\limsup_{n\to\infty}\|\mathfrak{W}_n\|_H^2\leqslant\mu_1\|\mathfrak{W}\|_{L^2(\Omega)}^2$

$\leqslant\|\mathfrak{W}\|_H^2,$

which together with the weak lower semi-continuity of the norm yields

$\limsup_{n\to\infty}\|\mathfrak{W}_n\|_H^2\leqslant\mu_1\|\mathfrak{W}\|_{L^2(\Omega)}^2\leqslant\|\mathfrak{W}\|_H^2\leqslant\liminf_{n\to\infty}\|\mathfrak{W}_n\|_H^2.$

This implies the strong convergence of $\{\mathfrak{W}_n\}$ in Hilbert space $H$, and hence $\|\mathfrak{W}\|_H=1$, and $\|\mathfrak{W}\|_H^2=\mu_1\|\mathfrak{W}\|_{L^2(\Omega)}^2$, which means that $\mathfrak{W}$ is the eigenfunction of the least positive eigenvalue $\mu_1$. Let $\varphi_1$ be the positive eigenfunction in the one-dimensional eigenfunction space of the least positive eigenvalue $\mu_1$ such that $\|\varphi_1\|_H=1$, and $\varphi_1(x)>0$ for all $x\in\Omega$, then $\mathfrak{W}=\pm\varphi_1$. Since $\mathfrak{W}_n(x)\to\mathfrak{W}(x)$, a.e. $x\in\Omega$, there exists $\Omega^*\subset\Omega$ such that $mes(\Omega/\Omega^*)=0$ and

$|\mathfrak{W}_n(x)|=\frac{|u_n(x)|}{\|u_n\|_H}\to\varphi_1(x),\quad\forall x\in\Omega^*,$

which implies $|u_n(x)|\to+\infty,\forall x\in\Omega^*$. Now, it follows by (3.43) that

$\mathfrak{P}\geqslant\int_{\Omega^*}\left[\frac{1}{2}\mu_1 u_n^2-F(u)\right]dx\to+\infty,\quad n\to\infty,$

which means that $\mathfrak{I}$ is bounded below.

By employing the similar methods used in [14], we shall prove that $\mathfrak{I}$ satisfies the (PS) condition.

In fact, if $\{u_n\}$ satisfies $\mathfrak{I}(u_n)\to a$, $\|\mathfrak{I}'(u_n)\|_{H^*}\to 0$, and $n$ is big enough, we get

$a+o(1)=\mathfrak{I}(u_n)=\frac{1}{2}\int_\Omega|\nabla u_n|^2 dx+\frac{1}{2}\frac{C}{D}\int_\Omega u_n^2 dx-\frac{A}{D}\int_\Omega F(u_n)dx$

and

$\langle\mathfrak{I}'(u_n),u_n\rangle=\int_\Omega|\nabla u_n|^2 dx+\frac{C}{D}\int_\Omega u_n^2 dx-\frac{A}{D}\int_\Omega f(u_n)u_n dx.$

And hence,

$\frac{1}{6}\|u_n\|_H^2+\frac{A}{D}\int_\Omega\left[\frac{1}{3}f(u_n)u_n-F(u_n)\right]dx$

$\leqslant a+o(1)+\frac{1}{3}\|\mathfrak{I}'(u_n)\|_{H^*}\|u_n\|_H,$

or

$\frac{1}{6}\|u_n\|_H^2\leqslant a+o(1)+\frac{1}{3}\|\mathfrak{I}'(u_n)\|_{H^*}\|u_n\|_H+\frac{A}{D}\int_\Omega|\frac{1}{3}f(u_n)u_n-F(u_n)|dx.$ (3.44)

On the other hand, (3.43) and (3.39) yield that there exists a constant $c_*>0$ big enough, satisfying

$\begin{cases}|\frac{1}{3}f(u)u-F(u)|\leqslant c_*u^{\frac{4}{3}}+c_*,&|u|\geqslant 1;\\|\frac{1}{3}f(u)u-F(u)|=|-\frac{1}{6}\frac{D}{A}\mu_1 u^2|\leqslant c_*,&u\in[-1,1],\end{cases}$

or

$|\frac{1}{3}f(u)u-F(u)|\leqslant c_*u^{\frac{4}{3}}+c_*,\quad\forall u\in\mathbb{R}^1,$

which together with (3.44), Holder inequality and Poincare inequality implies

$\frac{1}{6}\|u_n\|_H^2\leqslant a+o(1)+c_*\|u_n\|_H+\frac{c_*A}{D}mes(\Omega)+\frac{A}{D}c_*^2\|u_n\|_H^{\frac{4}{3}}.$ (3.45)

And then $\{u_n\}$ is bounded in $H$ due to (3.45).

Due to the fact $\|u\|^2\leqslant\|u\|_H^2\leqslant c_*\|u\|$, we see, $\|\cdot\|$ and $\|\cdot\|_H$ are a pair of equivalent norms. Moreover, we know from (3.39) that $f(u)$ satisfies the Caratheodory condition, and

$|f(u)|\leqslant c_*+c_*|u|^{\frac{1}{3}},\quad\text{and }0<\frac{1}{3}<2^*-1,$

which means that the bounded set $\{u_n\}$ with the condition $\|\mathfrak{I}'(u_n)\|_{H^*}\to 0$ is a compact set in the Hilbert space $H$. This have verified that $\mathfrak{I}$ satisfies the (PS) condition.

On the one hand, if $\inf_H\mathfrak{I}\geqslant 0$, we claim that there are infinitely many positive stationary solutions and infinitely many negative stationary solutions for the reaction diffusion system (3.43).

Indeed, since $\|\cdot\|$ and $\|\cdot\|_H$ are a pair of equivalent norms, Sobolev space $H$ has the orthogonal decomposition $H=E(\mu_1)\oplus E(\mu_1)^\perp$ [14,15], where $E(\mu_k)$ represents the eigenfunction space of $\mu_k$, and $E(\mu_1)^\perp=E(\mu_2)\oplus E(\mu_3)\oplus\cdots$. Obviously, $\mathfrak{I}$ satisfies (P1). In fact, if $u\in E(\mu_1)$ with $\|u\|_H\leqslant\delta$, the equivalence of norms in each finite dimensional space yields

$\|u\|_H\leqslant\delta\Rightarrow\int_\Omega|u(x)|dx<\delta_1\Rightarrow|u(x)|<1,\text{ a.e.}x\in\Omega.$

where the positive number $\delta$is small enough, and so is $\delta_1$. And then

$\mathfrak{I}(u)=\frac{1}{2}\|u\|_H^2-\frac{A}{D}\int_\Omega F(u)dx=\frac{1}{2}\|u\|_H^2-\mu_1\int_\Omega\frac{u^2}{2}dx=0$

$\leqslant 0,\quad u\in E(\mu_1),\|u\|_H\leqslant\delta.$ (3.46)

which together with $\inf_H\mathfrak{I}\geqslant 0$ implies that all $u\in E(\mu_1)$ with $\|u\|_H\leqslant\delta$ are the stationary solutions of the reaction diffusion system (3.38). Moreover, $E(\mu_1)$ implies that there are infinitely many positive stationary solutions and infinitely many negative stationary solutions for the reaction diffusion system (3.38).

On the other hand, if $\inf_H\mathfrak{I}<0$, we see, the condition (P3) holds. Then we claim that there are at least three equilibrium solutions including two non-zero stationary solutions, for the reaction diffusion system (3.38) if $\inf_H\mathfrak{I}<0$.

In fact, (3.46) implies that (P1) holds. In order to apply Lemma 2.2, we only need to verify the condition (P2) of Lemma 2.2. Next, for $u\in E(\mu_1)^\perp$, let $u=v+z$, where $v\in E(\mu_2),z\in E(\mu_3)\oplus E(\mu_4)\oplus\cdots$. Then we get

$\mathfrak{I}(u)\geqslant\left[\frac{1}{2}\mu_2\int_\Omega u^2 dx-\frac{A}{D}\int_\Omega F(u)dx\right]+\frac{1}{2}\left(1-\frac{\mu_2}{\mu_3}\right)\|z\|_H^2,\ \forall u\in E(\mu_1)^\perp.$ (3.47)

Due to (3.43), there exists $\delta\in(0,1)$ such that

$\frac{2A}{D}F(u)=\mu_1 u^2\leqslant\mu_2 u^2,\quad\text{if }|u|\leqslant\delta<1, u\in E(\mu_1)^\perp.$ (3.48)

Moreover, for this $\delta$, there exists correspondingly $\delta_2>0$ such that for $u\in E(\mu_2)$ with $\|u\|_H\leqslant\delta_2$, we get $|u(x)|\leqslant\frac{\delta}{2}$, a.e.$x\in\Omega$ in view of the equivalence of norms in finite dimensional space.

Define

$\Omega_1=\{x\in\Omega:|u(x)|\leqslant\delta\},\quad\Omega_2=\{x\in\Omega:|u(x)|>\delta\}.$

Due to the orthogonal decomposition of the Sobolev space $H$ and $u=v+z$, we see, $\|u\|_H\leqslant\delta_2\Rightarrow\|v\|_H\leqslant\delta_2$, which implies

$|v(x)|\leqslant\frac{\delta}{2}\leqslant\frac{1}{2}|u(x)|,\quad|z(x)|\geqslant|u(x)|-|v(x)|\geqslant\frac{1}{2}|u(x)|,\text{ a.e.}x\in\Omega_2.$

Besides, (3.43) yields

$|\frac{1}{2}\mu_2 u^2 dx-\frac{A}{D}F(u)|\leqslant c_*|u|^3\leqslant 8c_*|z|^3, x\in\Omega_2.$

So we can see it from the orthogonal decomposition of the Sobolev space $H$ and the Sobolev embedding theorem that for all $u \in E(\mu_1)^{\perp}$, if $\|u\|_H \leqslant \delta_2$, we get

$$\Im(u) \geqslant \frac{1}{2}\left(1 - \frac{\mu_2}{\mu_3}\right)\|z\|_H^2 - 8c_*^2\|z\|_H^3,$$

which implies that the condition (P2) holds. And now all the conditions of Saddle point theorem are satisfied. According to Lemma 2.2, we have proved the claim. And the proof is completed.  □

**Remark 11.** Statement 2 actually illuminates that under Lipschtiz assumption on active function, diffusion may make an equilibrium point become infinitely many stationary solution. The deeper purpose of Statement 2 may be revealed in the final section for further considerations (see Problem 1). On the other hand, according to the Introduction in [20], the function $f$ defined by (3.39) satisfies the

$$
\begin{cases}
\frac{\partial y(t,x)}{\partial t} = D_1\Delta y(t,x) - C_1 y(t,x) + A_1 g(y(t,x)) + B_1 g(y(t - \tau(t),x)) + J_1, & (t,x) \in \mathbb{R}_+ \times \Omega_1, \\
y_i(t,x) = 0, t \geqslant 0, x \in \partial\Omega_1, \ i = 1,2,
\end{cases}
\tag{4.1a}
$$

$$
\begin{cases}
\frac{\partial y(t,x)}{\partial t} = D_2\Delta y(t,x) - C_2 y(t,x) + A_2 g(y(t,x)) + B_2 g(y(t - \tau(t),x)) + J_2, & (t,x) \in \mathbb{R}_+ \times \Omega_2, \\
y_i(t,x) = 0, t \geqslant 0, x \in \partial\Omega_2, \ i = 1,2,
\end{cases}
\tag{4.1b}
$$

and

$$
\begin{cases}
\frac{\partial y(t,x)}{\partial t} = D_3\Delta y(t,x) - C_3 y(t,x) + A_3 g(y(t,x)) + B_3 g(y(t - \tau(t),x)) + J_3, & (t,x) \in \mathbb{R}_+ \times \Omega_3, \\
y_i(t,x) = 0, t \geqslant 0, x \in \partial\Omega_3, \ i = 1,2,
\end{cases}
\tag{4.1c}
$$

or the following corresponding homogeneous equations:

conditions [20,7,20,8], and hence, the zero solution is actually the unique equilibrium point of the ordinary differential system (3.37). Thereby, Statement 2 has actually verified that under the Lipschitz assumption on the function $f$, the small diffusion can truly make one equilibrium point become multiple stationary solutions (three stationary solutions, even infinitely many stationary solutions).

Finally, in this section, if I want to obtain the existence of exponentially stable common stationary solution rather than a positive stationary solution of the system (2.1), I only need to replace the condition (A2) with the following condition:

(A2*) There is a positive real number $c > 0$ such that

$$| -C_\sigma v + A_\sigma g(v) + B_\sigma g(v) + J_\sigma | \leqslant cD_\sigma E, \quad \forall v \in R^n$$

where $D_\sigma > 0$ is a positive definite diagonal matrix, and $E = (1, 1, \cdots, 1)^T \in \mathbb{R}^n$.

In fact, I only need to replace the definition of $\Re$ in the proof of Theorem 3.1 with

$$\Re = \left\{\varphi(x) \in \left[C\left(\overline{\Omega_\sigma}\right)\right]^n : \|\varphi(x)\| < +\infty; \varphi(x) = 0, x \in \partial\Omega\right\}, \tag{3.49}$$

which implies that $\Re$ is a closed convex set. Then Theorem 3.1 and its proof have got the following direct corollary:

**Corollary 3.6.** *Suppose that the conditions (A1) and (A2*) hold, then the system (2.1) possesses a bounded stationary solution $y^\sigma(x)$ for $x \in \Omega_\sigma$ with $y^\sigma|_{\partial\Omega_\sigma} = 0$. If, in addition, the condition (3.1) holds, then the null solution of the switched delayed reaction–diffusion system (2.3) equipped with the initial value (2.4) is exponentially stable with the convergence rate $\frac{\gamma}{2}$.*

It is natural that Theorem 3.2 has got the following direct corollary:

**Corollary 3.7.** *If all the assumptions of Corollary 3.6 hold, and if, in addition, the condition (A3) is satisfied, then the system (2.1) possesses a unique bounded stationary solution $y^\sigma(x)$ for $x \in \Omega_\sigma$ with $y^\sigma|_{\partial\Omega_\sigma} = 0$. And the null solution of the switched delayed reaction–diffusion system (2.3) equipped with the initial value (2.4) is globally exponentially stable with the convergence rate $\frac{\gamma}{2}$. Particularly, if $\mathfrak{I} = \{1\}$ or $N = 1$, the unique stationary solution $y^\sigma(x)(\sigma = 1)$ of the **deterministic system** (2.1) is globally exponentially stable **in the classical meaning**.*

## 4. Numerical example

**Example 4.1.** Consider the following switched financial system with $\sigma \in \{1, 2, 3\}$,

equipped with the initial value:

$$u_j(s,x) = \phi_j(s,x) = \prod_{\sigma=1}^3 \sin^j\left[x_1^{33}(x_1 - 5(\sigma+1))^{353}x_2^{63}(x_2 - 5(\sigma+1))^{79}\right],$$

$$j = 1,2, -\tau \leqslant s \leqslant 0, x \in \Omega_\sigma, \ \sigma \in \{1,2,3\},$$
$$\tag{4.3}$$

where $\Omega_1 = [0,1] \times [0,1], \Omega_2 = [0,1.3] \times [0,1.3], \Omega_3 = [0,1.5] \times [0,1.5], \lambda_{11} = 19.7392, \lambda_{21} = 11.68, \lambda_{31} = 8.7730$ (see Remark 14).

Set $c = 100000000, I = diag(1,1)$, and

$$C_1 = \begin{pmatrix} 0.448 & 0 \\ 0 & 0.441 \end{pmatrix}, C_2 = \begin{pmatrix} 0.455 & 0 \\ 0 & 0.441 \end{pmatrix},$$

$$C_3 = \begin{pmatrix} 0.438 & 0 \\ 0 & 0.433 \end{pmatrix},$$

$$A_1 = \begin{pmatrix} 0.45 & 0.00003 \\ -0.00003 & 0.44 \end{pmatrix}, A_2 = \begin{pmatrix} 0.452 & 0.00001 \\ -0.00001 & 0.441 \end{pmatrix},$$
$$A_3 = \begin{pmatrix} 0.439 & 0.000015 \\ -0.00001 & 0.433 \end{pmatrix},$$

$$B_1 = \begin{pmatrix} 0.446 & -0.00003 \\ 0.00003 & 0.442 \end{pmatrix}, B_2 = \begin{pmatrix} 0.458 & -0.00001 \\ 0.00001 & 0.441 \end{pmatrix},$$
$$B_3 = \begin{pmatrix} 0.437 & -0.000015 \\ 0.00001 & 0.433 \end{pmatrix},$$

Let $g_i(y_i) = \frac{39 + 2y_i + 0.000001\sin y_i}{4}$, and then $G = diag(0.51, 0.51)$. Set $J_\sigma = (0.2\sin\sigma, -0.1\cos\sigma)^T, \sigma = 1,2,3$, then the direct calculation can verify that both conditions (A1) and (A2) hold.

For example, in the case of $\sigma = 1$, direct computation derives

$$0 \leqslant [-C_1 v + A_1 g(v) + B_1 g(v) + J_1] = \begin{pmatrix} 0.448 \times \frac{39 + 0.000001 \sin v_1}{2} \\ 0.441 \times \frac{39 + 0.000001 \sin v_2}{2} \end{pmatrix}$$

$$+ \begin{pmatrix} 0.2 \sin 1 \\ -0.1 \cos 1 \end{pmatrix} \leqslant c D_1 E,$$

which means that (A2) holds in the case of $\sigma = 1$. Similarly we can compute and verify that (A2) holds in the case of $\sigma = 2, 3$ (see Figs. 1 and 2).

**Case 1**

Set

$$D_1 = \begin{pmatrix} 0.05 & 0 \\ 0 & 0.055 \end{pmatrix}, D_2 = \begin{pmatrix} 0.07 & 0 \\ 0 & 0.075 \end{pmatrix}, D_3 = \begin{pmatrix} 0.09 & 0 \\ 0 & 0.095 \end{pmatrix}, \tag{4.1d}$$

Moreover, in (3.17), let $p_\sigma \equiv 1$, then it is easy to verify that $p_\sigma > 0$ satisfies $p_\sigma^2 I \geqslant (A_\sigma + B_\sigma)^T (A_\sigma + B_\sigma)$ for each $\sigma$. In addition, set $\varepsilon = 2$,

$$- C_\sigma + \frac{p_\sigma}{2} \left( \varepsilon^{-1} I + \varepsilon G^2 \right) < 0.1 I < \lambda_{\sigma 1} D_\sigma, \quad \text{for all } \sigma$$

Then the condition (A3) holds (Below, it can be verified similarly that (A3) holds in Case 2–3, too).

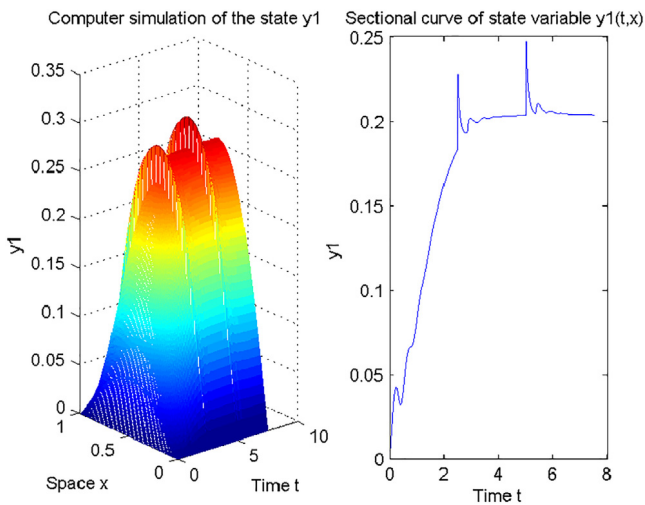

**Fig. 1.** Computer simulation of the state variable y1

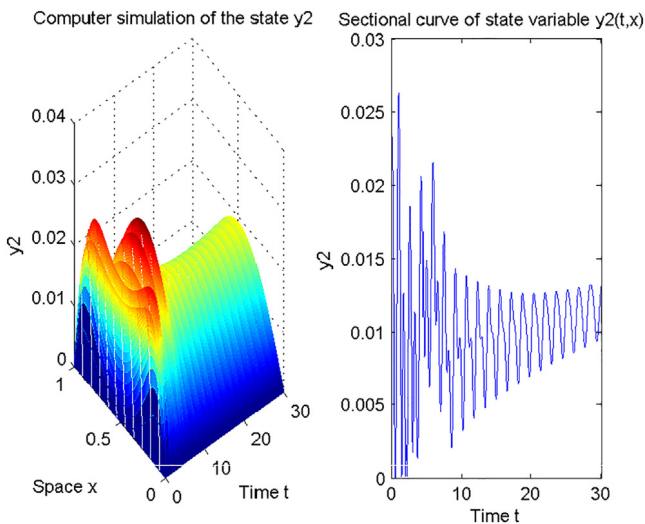

**Fig. 2.** Computer simulation of the state variable y2

Let $q = 1.00001, \Psi = 0.00018 I$ and $\tau = 3.5$, then employing computer LMI toolbox to solve the inequality (3.1) derives the following feasible data:

$$\beta_1 = 0.5676, \ \beta_2 = 0.3633, \ \beta_3 = 0.0691, \ \gamma = 0.38,$$

then the switched delayed reaction–diffusion system (2.3) equipped with the initial value (2.4) is globally exponentially stable with the convergence rate 19% due to Theorem 3.2.

**Case 2**

If replacing the diffusion coefficients (4.4) with the following diffusion coefficients

$$D_1 = \begin{pmatrix} 0.1 & 0 \\ 0 & 0.15 \end{pmatrix}, D_2 = \begin{pmatrix} 0.15 & 0 \\ 0 & 0.2 \end{pmatrix}, D_3 = \begin{pmatrix} 0.1 & 0 \\ 0 & 0.15 \end{pmatrix}, \tag{4.1e}$$

and other data of Case 1 are not changed, we can use the computer LMI toolbox to solve the inequality (3.1), resulting in the following feasible data:

$$\beta_1 = 0.6769, \ \beta_2 = 0.2333, \ \beta_3 = 0.0898, \ \gamma = 0.44,$$

then the switched delayed reaction–diffusion system (2.3) equipped with the initial value (2.4) is globally exponentially stable with the convergence rate 22% due to Theorem 3.2.

**Case 3**

If replacing $\tau = 3.5$ with $\tau = 3$, and other data of Case 1 are not changed, we can use the computer LMI toolbox to solve the inequality (3.1), resulting in the following feasible data:

$$\beta_1 = 0.6616, \ \beta_2 = 0.3113, \ \beta_3 = 0.0271, \ \gamma = 0.58,$$

then the switched delayed reaction–diffusion system (2.3) equipped with the initial value (2.4) is globally exponentially stable with the convergence rate 29% due to Theorem 3.2.

**Remark 12.** Table 1 tells us that the larger the diffusion coefficient, the faster the convergence rate. On the other hand, the harsh condition (A2) illuminates that the diffusion makes it more difficult to judge the stability of the system.

**Remark 13.** Table 2 indicates that the larger the upper bound $\tau$ of time delays, the slower the convergence speed $\frac{\gamma}{2}$.

Example 4.1 shows the feasibility of our [Theorem 3.1–3.2], whose unique existence criteria of positive stationary solution were successfully extended to ecosystem ([19, Theorem 3.1–3.2]).

**Remark 14.** The smallest positive eigenvalue of $-\Delta_p$ in $W_0^{1,p}(0, T)$ is

**Table 1**

Comparisons the influences on the convergence rate $\frac{\gamma}{2}$ under different diffusion coefficients with the same other data.

| | Case 1 | Case 2 |
|---|---|---|
| Diffusion coefficient | (4.4) (smaller) | (4.5) (bigger) |
| Convergence rate | 19% | 22% |

**Table 2**

Comparisons the influences on the convergence rate $\frac{\gamma}{2}$ under different upper limits of delays with the same other data.

| | Case 1 | Case 3 |
|---|---|---|
| Diffusion coefficient | $\tau = 3.5$ | $\tau = 3$ |
| Convergence rate | 19% | 29% |

$$\lambda_1 = \left( \frac{2}{T} \int_0^{(p-1)^{\frac{1}{p}}} \frac{dt}{\left(1 - \frac{t^p}{p-1}\right)^{\frac{1}{p}}} \right)^p.$$

If $\Omega = \left\{ (x_1, x_2)^T : 0 < x_1 < \alpha, 0 < x_2 < \beta \right\} \subset \mathbb{R}^2$ and $W_0^{1,p}(\Omega)$ with $p = 2$, the first eigenvalue $\lambda_1 = \left(\frac{\pi}{\alpha}\right)^2 + \left(\frac{\pi}{\beta}\right)^2$.

Besides, it is well known that there is the following approximate substitution of Poincare inequality lemma:

**Remark 15.** Let $\Omega$ be a cube $|x_i| < l_i (i = 1, 2, \cdots, n)$ and let $\mu(x)$ be a real-valued function belonging to $C^1(\Omega)$ which vanish on the boundary $\partial \Omega$ of $\Omega$, i.e., $\mu(x)|_{\partial \Omega} = 0$, then

$$\int_\Omega \mu^2(x) dx \leqslant l_i^2 \int_\Omega |\frac{\partial \mu}{\partial x_i}|^2 dx.$$

## 5. Conclusions and further considerations

By constructing a compact operator on a convex set, the author makes up for the loss of compactness in infinite dimensional space. Using a fixed point theorem, variational methods and Lyapunov functional method results in the existence positive bounded stationary solution, which is exponentially stable. Moreover, by using the first positive eigenvalue of Laplace operator $-\Delta$ to restrain Lipschitz constants, the author proposes the uniqueness theorem of the stationary solution of reaction–diffusion system under Dirichlet zero boundary value, and thereby the stability of Theorem 3.1 becomes global. Not only that, Theorem 3.1 and Theorem 3.2 derive a corollary on the variance of global stability. Moreover, Statement 1 points out the fact that non-zero constant vector can not be a stationary solution of the reaction–diffusion system under Dirichlet zero boundary value. Besides, Statement 2 points out that the influence of diffusions changes the number of the system under Lipschitz assumptions on active functions. Finally, a numerical example is presented to illuminate the effectiveness of the proposed methods.

Below, some interesting problems are proposed as follows.

**Problem 1.** How to improve the example in Statement 2 or add a suitable condition to the example in Statement 2 so that the constant equilibrium point of the ordinary differential system can be truly proved to be globally asymptotically stable. At the same time, the corresponding reaction–diffusion system owns multiple stationary solutions. If so, many global stability results of delayed neural networks in the form of ordinary differential equations may only be locally asymptotical stability criteria in real engineering due to the inevitable diffusions.

**Problem 2.** How to replace (A2) with a weaker condition in Theorem 3.1?

**Problem 3.** Is the condition (A3) of Theorem 3.2 necessary for the uniqueness? If not, what's the weaker condition?

**Problem 4.** In Corollary 3.4 (Global Stability Invariance), $u^*$ of Conclusion (1) becomes $u^*(x)$ of Conclusion (2) due to the small diffusion (the diffusion coefficients $\{r_i\}$). Now I wonder whether for any given $\varepsilon > 0$, there exits a corresponding number $\delta_\varepsilon > 0$ such that $\|u^*(x) - u^*\| < \varepsilon$ or $\sup\limits_{x \in \Omega} |u^*(x) - u^*| < \varepsilon$ if $\max\limits_i |r_i| < \delta_\varepsilon$. On the other hand, $u^*(x)$ is dependent on $\{r_i\}$. Now I wonder whether $u^*(x)$ is a continuous dependence on $\max\limits_i |r_i|$ or $\{r_i\}$.

**Problem 5.** How to give a boundedness criterion about the state variable of the neural networks (2.1) similarly as that of [46, Theorem 3.3], to some extent?

**Remark 16.** There are a lot of boundedness results for reaction–diffusion neural networks in previous literature (see, e.g. [6,29,39] and their references therein). For example, in [29, Theorem 1], the boundedness of nonautonomous reaction–diffusion neural networks was obtained in the meaning of the norm $\|\cdot\|_p$ with $p \neq \infty$. However, about Problem 5, in view of the viewpoint in Remark 1 of this paper, we need the boundedness of the state variables in the meaning of the norm $\|\cdot\|_\infty$, similarly as that of [46, Theorem 3.3] to some extent.

Recently, Xiaodi Li and his coauthors in [41,42] proposed some methods on dealing with the infinite delay, deducing a sequence of interesting results. Finally in this paper we propose the following problem:

**Problem 6.** Reconsidering all the results obtained in this paper under the infinite delays, which might means a series of new results. This is an interesting problem, too.

### Funding Statement

The work is supported by the Application basic research project of science and Technology Department of Sichuan Province (No. 2020YJ0434) and the Major scientific research projects of Chengdu Normal University in 2019 (No. CS19ZDZ01)

### Declaration of Competing Interest

The authors declare that they have no known competing financial interests or personal relationships that could have appeared to influence the work reported in this paper.

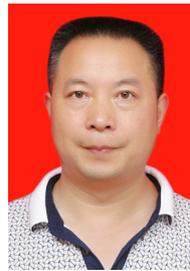

**Ruofeng Rao**, born in 1969, male, professor. From 2003 to 2007, he is a college lecturer in Jiujiang University, Zhaoqing University, and Yibin University, China. From 2008 to 2013, he is an associate professor at Yibin University, and since the end of 2013, he is a professor at Yibin University. And he is currently a professor at Chengdu Normal University. He is a Guest Editor of a international academic journal, and a reviewer of more than 8 other journals. He is the author or coauthor of more than 80 journal papers, published in international academic journals or Chinese academic journals. His research interests include variation methods, critical point theory, stochastic control, fixed point theory and its applications.

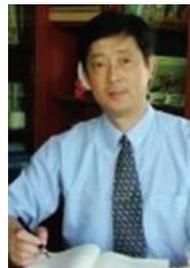

**Jialin Huang**, born in 1949, male, professor, his research interests include mathematics and applied mathematics. Now he teaches in the basic Department of Sichuan Sanhe vocational college.

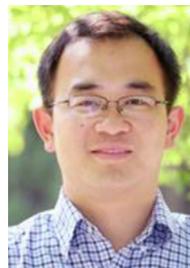

**Xiaodi Li** received the B.S. and M.S. degrees from Shandong Normal University, Jinan, China, in 2005 and 2008, respectively, and the Ph.D. degree from Xiamen University, Xiamen, China, in 2011, all in applied mathematics. He is currently a Professor with the School of Mathematics and Statistics, Shandong Normal University. From Nov. 2014 to Dec. 2016, he was a Visiting Research Fellow at Laboratory for Industrial and Applied Mathematics in York University, Canada, and the University of Texas at Dallas, USA. In 2017 and 2019, he was working as Visiting Research Fellow at the Department of Mathematics, City University of Hong Kong, and Department of Mechanical Engineering, The University of Hong Kong, Hong Kong. respectively. He has authored or coauthored more than 70 research papers. He is currently the Editor-in-Chief of the journal 'AIMS Mathematics' and an editor of 'AAM', 'IJDE' etc. His research interests include control theory, hybrid systems, time-delay systems, neural networks, and applied mathematics.